\documentclass[preprint,12pt]{elsarticle}

\usepackage{amssymb}
\usepackage{amsmath}
\usepackage{amsfonts}
\usepackage{amsopn}
\usepackage{hyperref}
\usepackage{color}

\usepackage{cancel}
\usepackage{pifont}
\usepackage{multirow}
\usepackage{bm}
\usepackage{booktabs}
\usepackage[normalem]{ulem} 
\graphicspath{{./content/figures}}

\newcommand{\revisea}[1]{{#1}} 
\newcommand{\reviseb}[1]{{#1}} 
\newcommand{\revisec}[1]{{#1}} 
\newcommand{\appref}[1]{#1}

\usepackage{times}

\journal{Physica D}

\begin{document}

\begin{frontmatter}

\title{Offline Supervised Learning V.S. Online Direct Policy Optimization: A Comparative Study and A Unified Training Paradigm for Neural Network-Based Optimal Feedback Control}

\author[label1]{Yue Zhao}
\ead{me@y-zh.com}
\author[label2]{Jiequn Han}
\ead{jiequnhan@gmail.com}

\affiliation[label1]{organization={Center for Data Science, Peking University},%Department and Organization
            addressline={No. 5 Yiheyuan Road}, 
            city={Beijing},
            postcode={100871}, 
            % state={},
            country={China}}

\affiliation[label2]{organization={Flatiron Institute},%Department and Organization
            addressline={162 5th Ave.}, 
            city={New York},
            postcode={10010}, 
            state={NY},
            country={USA}}

\begin{abstract}
This work is concerned with solving neural network-based feedback controllers efficiently for optimal control problems. We first conduct a comparative study of two prevalent approaches: \textit{offline supervised learning} and \textit{online direct policy optimization}. Albeit the training part of the supervised learning approach is relatively easy, the success of the method heavily depends on the optimal control dataset generated by open-loop optimal control solvers.
In contrast, direct policy optimization turns the optimal control problem into an optimization problem directly without any requirement of pre-computing, but the dynamics-related objective can be hard to optimize when the problem is complicated.
Our results underscore the superiority of offline supervised learning in terms of both optimality and training time. To overcome the main challenges, \textit{dataset} and \textit{optimization}, in the two approaches respectively, we complement them and propose the \textbf{Pre-train and Fine-tune} strategy as a unified training paradigm for optimal feedback control, which further improves the performance and robustness significantly. Our code is accessible at \url{https://github.com/yzhao98/DeepOptimalControl}.
\end{abstract}

% %%Research highlights
% \begin{highlights}
%  \item We establish a benchmark and comprehensively compare \textit{offline supervised learning} and \textit{online direct policy optimization} for solving neural network-based optimal feedback control problems. 
%  \item We identify the challenges of two methods, which are dataset and optimization, respectively.
% \item Our results underscore the superiority of offline supervised learning to online direct policy optimization in terms of both optimality and training time.
% \item Drawing from our comparative analysis, we propose a paradigm to train neural network-based feedback controllers, namely \textbf{Pre-Train} \textbf{and} \textbf{Fine-Tune}, which significantly enhances performance and robustness.
% \end{highlights}


\end{frontmatter}

\section{Introduction}

It is ubiquitous and paramount to design optimal feedback controllers~\citep{franklinFeedbackControlDynamic2020} for various complex tasks in engineering and industry. Real-world applications are more than challenging due to high-dimensionality and speed requirements for real-time execution. In recent years, deep learning has been introduced to tackle these issues and shown astonishing performances~\citep{hanDeepLearningApproximation2016, Adaptive_BVP_HJB, AIPontryagin, Empowering}. Generally speaking, given system dynamics explicitly, there are two prevalent types of approaches to training neural network-based optimal feedback controllers: \textit{offline supervised learning} and \textit{online direct policy optimization}. 

Offline supervised learning is to train a feedback network controller by approximating the corresponding open-loop solutions at different states directly, leveraging the fact that open-loop control for a fixed initial state is much easier to solve than feedback control in optimal control problems (OCP). 

The other way, online direct policy optimization, 
transforms the network-based feedback OCP concerning a distribution of initial states into an optimization problem and solves it directly without pre-computing open-loop optimal control. We call this method direct policy optimization because the objective to minimize comes exactly from the original OCP.  
With fully-known dynamics, the evolution of states (controlled by a neural network) is governed by a controlled Ordinary Differential Equation (ODE), and one can use stochastic gradient descent to optimize network parameters. This type of approach was first proposed to solve high-dimensional stochastic control problems~\citep{hanDeepLearningApproximation2016} and recently applied to solve deterministic optimal control problems~\citep{AIPontryagin, CTPG}.  

The two methods have not, though, hitherto been comparatively studied.
Thus, it is imperative to conduct a comparative study of them. 
Albeit supervised learning problems can be easily optimized, the quality of the learned controller heavily depends on the dataset generated by open-loop optimal control solvers.
In this work, we demonstrate that the supervised learning approach holds advantages over direct policy optimization in terms of optimality and training time. The primary benefit of supervised learning stems from the network's ability to directly learn the optimal control signal from a pre-computed dataset, making the process less challenging compared to the online training required for direct policy optimization. However, direct policy optimization can yield further performance enhancements with a near-optimally initialized network, since its objective takes the controlled dynamics into consideration and aligns more closely with the original goal of the OCP.

Our comparison of the two approaches shows that they are interrelated and complementary, enlightening us to combine them to overcome their inherent limitations. Thereby, we present a unified training paradigm for neural network-based optimal feedback controllers, dubbed as the \textbf{Pre-Train and Fine-Tune} strategy. Initial pre-training through \textit{offline supervised learning} guides the network to a reasonable solution with a small loss, whilst fine-tuning by \textit{online direct policy optimization} breaks the limitations of the precomputed dataset and further improves the controller's performance and robustness.

We summarize our contributions as follows:

\begin{itemize}
 \item We establish a benchmark and comprehensively compare \textit{offline supervised learning} and \textit{online direct policy optimization} for solving neural network-based optimal feedback control problems. 
 \item We identify the challenges of two methods, which are \revisea{dataset quality and network optimization}, respectively.
\item Our results underscore the superiority of offline supervised learning to online direct policy optimization in terms of both optimality and training time.
\item Drawing from our comparative analysis, we propose a \revisea{\sout{new}} paradigm to train neural network-based feedback controllers, namely \textbf{Pre-Train} \textbf{and} \textbf{Fine-Tune}, which significantly enhances performance and robustness. 
\end{itemize}

\section{Preliminaries and Related Works}

\subsection{Mathematical Formulation}

We first consider an open-loop optimal control problem:
\begin{equation}
\left\{\begin{array}{cl}
\underset{\boldsymbol{u}(t)}{\operatorname{minimize}} & J[\boldsymbol{u}(t)]=M\left(\boldsymbol{x}\left(T\right)\right)+\displaystyle{\int_0^{T}} {L}(\boldsymbol{x}, \boldsymbol{u}) d t, \\
\text { s.t. } & \dot{\boldsymbol{x}}(t)=\boldsymbol{f}(\boldsymbol{x}, \boldsymbol{u}), \quad \boldsymbol{x}(0)=\boldsymbol{x}_0,\\
% & 
\end{array}\right .
\label{eq1-objective}
\end{equation}
where $\boldsymbol{x}_0\in \mathbb{R}^n$ denotes an initial state, $\boldsymbol{x}(t) \in \mathbb{R}^n$ denotes the state at time $t \in [0, T]$, $\mathcal{U} $ denotes the admissible control set, and $\boldsymbol{u}(t) \in \mathcal{U} \subset \mathbb{R}^m$ denotes the open-loop control function. The dynamics is described by a smooth function $\boldsymbol{f}: \mathbb{R}^n \times \mathcal{U} \rightarrow \mathbb{R}^n$. The total cost $J$ is the sum of the terminal cost $M: \mathbb{R}^n \rightarrow \mathbb{R}$ and an integral of the running cost $L:\mathbb{R}^n \times \mathcal{U} \rightarrow \mathbb{R}$, \revisea{with the assumption that both $M$ and $L$ are differentiable}.\footnote{For ease of notation, we assume the dynamics and running cost are time-independent. Both approaches discussed in this paper can be applied straightforwardly to problems with time-dependent dynamics or running cost.} We assume the solution to the open-loop OCP~(\ref{eq1-objective}) exists and is unique, which is  $\boldsymbol{u}^{*}:[0, T] \rightarrow \mathcal{U}$ that minimizes the total cost under the dynamics and the given initial state.

The open-loop control is designed for a fixed initial state and takes only time $t$ as the input. In contrast, the closed-loop control is designed for various initial states and thus takes both time and the current state as the input, \textit{i.e.}, $\boldsymbol{u}(t, \boldsymbol{x}):[0, T] \times \mathbb{R}^n \rightarrow \mathcal{U}\subset \mathbb{R}^m$. In classical control theory~\citep{Openloop}, it is well-known that there exists a closed-loop optimal control $\boldsymbol{u}^{*}:[0, T]\times \mathbb{R}^n \rightarrow \mathcal{U}$ so that $\boldsymbol{u}^*\left(t, \boldsymbol{x}^*\left(t ; \boldsymbol{x}_0\right)\right)$ is identical to the open-loop optimal control solution at time $t$ with the initial state $\boldsymbol{x}_0$, where $\boldsymbol{x}^*$ follows $\dot{\boldsymbol{x}}^*\left(t ; \boldsymbol{x}_0\right)=\boldsymbol{f}\left(\boldsymbol{x}^*\left(t ;\boldsymbol{x}_0\right), \boldsymbol{u}^*\left(t, \boldsymbol{x}^{*}\left(t ; \boldsymbol{x}_0\right)\right)\right.$. The identity holds for any initial state $\boldsymbol{x}_0$. In other words, we can induce a family of open-loop optimal control with any possible initial states if the closed-loop optimal control is given. Based on that fact, we slightly abuse the notation $\boldsymbol{u}$ to denote both the closed-loop control and the induced open-loop control. The context of closed-loop or open-loop control can be inferred from the arguments when necessary and will not be confusing. We also use the words feedback and closed-loop interchangeably. We focus on the design of closed-loop optimal controllers,
considering they are more reliable and robust to dynamic disturbance and model misspecification in real-world applications compared to the open-loop counterpart. \revisea{We operate under the assumption that there is a reliable numerical method for deriving the open-loop optimal solution, which implies that the problem is not overly complex, thus making the development of closed-loop optimal control feasible.}

Traditional methods for solving the closed-loop optimal control rely on solving the corresponding Hamilton--Jacobi--Bellman (HJB) equation. It is notoriously difficult to solve those equations in high dimensions based on classical grid-based methods, due to the so-called {\it{curse of dimensionality}}~\citep{bellmanDynamicProgramming1957}. \revisec{To overcome this essential difficulty, there has been active research in recent years to approximate the control and value functions with other function approximators and proper objective functions to find the optimal approximators. Notable examples of such approximating functions include sparse polynomials~\cite{azmi2021optimal, kunisch2023learning}, sparse grids~\cite{Satellite_17}, kernel functions~\cite{weston2002kernel}, and neural networks~\cite{hanDeepLearningApproximation2016, Adaptive_BVP_HJB, AIPontryagin}. Particularly, neural networks have received significant interest due to their exceptional ability in high-dimensional function approximation and their flexibility in optimization with various loss functions. For instance, it is feasible to employ loss functions linked to the classical Hamilton-Jacobi-Bellman (HJB) equation~\cite{kunisch2023learning} or the Pontryagin principle~\cite{meng2022sympocnet,onken2022neural}, especially when the Hamiltonian minimization can be explicitly solved. In this work, we focus on neural network-based methods and loss functions that do not require the previously mentioned condition. The progress made in this direction can be categorized into two families of methods: \textit{offline supervised learning} and \textit{online direct policy optimization}. 
We introduce their details in the following subsections accordingly.}

\subsection{Offline Supervised Learning}

In the supervised learning approach, neural networks are utilized to approximate solutions provided by an optimal control dataset. The learning signal \textit{de facto} comes from the open-loop optimal control solution, which is much easier to solve than the closed-loop control. It is offline trained and then applied in real-time control. The regression problem is easily formulated and optimized. The key point of this training approach is how to generate high-quality data for training the neural network to achieve optimal control.

\paragraph{Data Generation} Many numerical methods can be applied to solve the open-loop optimal control, such as numerical solvers for the Boundary Value Problem (BVP,~\citealp{BVP}) and Differential Dynamical Programming (DDP,~\citealp{DDP}).  
Since the data is generated offline, the computation time of each trajectory is usually not the primary concern, as long as it is affordable to generate a suitable dataset in the form of
$
\mathcal{D} = \{ (t^{(i)}, \boldsymbol{x}^{(i)}, \boldsymbol{u}^{(i)} ) \},
$
where $\boldsymbol{u}^{(i)}$ denotes the optimal control at time $t^{(i)}$ and state $\boldsymbol{x}^{(i)}$, with $i$ being the index of data. 
For instance,~\cite{Adaptive_BVP_HJB} solves the corresponding BVP as a necessary condition of the optimal control based on PMP and collects open-loop solutions starting from different initial states and evaluated at different temporal grid points to form the dataset.

\paragraph{Objective} The objective in offline supervised learning is to find a network policy $\boldsymbol{u}^{\text{NN}}(t, \boldsymbol{x}; \boldsymbol{\theta})$ with parameter $\boldsymbol{\theta}$ that minimizes the least-square error based on the generated dataset $\mathcal{D}$:

\begin{equation}
\underset{\boldsymbol{\theta}}{\operatorname{minimize}}\quad \frac{1}{|\mathcal{D}|} \sum_{i=1}^{|\mathcal{D}|}\left\|\boldsymbol{u}^{\mathrm{NN}}\left(t^{(i)},\boldsymbol{x}^{(i)} ; \boldsymbol{\theta}\right) - \boldsymbol{u}^{(i)} \right\|^2.
\end{equation}

\revisec{We remark that this approach of learning the feedback control is similar to imitation learning~\cite{osa2018algorithmic} or behavior cloning~\cite{bain1995framework}, where the objective is to learn from a dataset with labels from expert policies. Our method differs from conventional imitation learning by how to obtain expert demonstrations. Instead of human demonstrations, labels are derived from solving corresponding open-loop optimal control problems, enabling the construction of an expansive dataset with optimal quality that can be scaled up. Also, this approach is not limited to learning control policies but is also applicable to value functions~\cite{Adaptive_BVP_HJB}.
}

\subsection{Online Direct Policy Optimization}\label{sec: direct}

The direct method~\citep{Bock_Direct, Betts_Direct, Ross_Direct, Diehl_Direct} is a classical family of methods in the optimal control literature that mainly refers to methods that transform the open-loop OCP into a nonlinear optimization problem. The neural network-based direct policy optimization investigated in this paper shares a similar optimization objective with those direct methods, but distinguishes from them in terms of solving closed-loop optimal controls. That is why we highlight the term ``policy optimization" in its name. The method applies to both stochastic and deterministic control problems, and in this work we focus on problems with deterministic dynamics and cost.

\paragraph{Objective} 
To formulate the direct policy optimization method, we assume $\boldsymbol{x}_0$ is sampled from a distribution ${\mu}$ that covers the initial state space of interest, such as a uniform distribution over a compact set. We use $J(\boldsymbol{u}; \boldsymbol{x_0})$ to denote the cost defined in \eqref{eq1-objective}, which allows the control to be feedback and highlights the effect of the initial state. The objective function in the direct policy optimization is defined as an expectation of the cost functional over the distribution of initial states: 

\begin{equation}
\underset{\boldsymbol{\theta}}{\operatorname{minimize}}\quad \underset{\boldsymbol{x}_0 \sim {\mu}}{\mathbb{E}} J(\boldsymbol{u}^{\mathrm{NN}}(\cdot; \boldsymbol{\theta}); \boldsymbol{x}_0).
\label{eq3:obj-direct}
\end{equation}
We classify direct policy optimization into two settings based on whether the system is explicitly known: the fully-known dynamics setting and the unknown dynamics setting.

With fully-known dynamics, we have explicit forms of $\eqref{eq1-objective}$ (the dynamics $f$, cost functions $M$ and $L$) and build a large computational graph along trajectories corresponding to the objective. Then we can optimize the objective based on randomly sampled initial states to approximate the expectation. \revisec{Employing fixed initial points is also a viable approach, as in~\cite{kunisch2023learning, kunisch2021semiglobal}; however, in complex scenarios, randomly sampled initial points better cover the initial distribution, leading to improved stability and performance.} The gradient with respect to the network parameters $\boldsymbol{\theta}$ can be computed using two methods: back-propagation~\citep{bp} or the adjoint method~\citep{CTPG, pontryaginMathematicalTheoryOptimal1962}. Details of these methods and a comparison between them are provided in \appref{\ref{app:adjoint}}. We report results optimized using back-propagation in the main content, as the adjoint method is much more time-consuming with indistinguishable objective improvements.

In contrast, in the setting with unknown dynamics, we can only observe trajectories generated by the underlying dynamics and cost signals on the encountered states. This scenario is commonly referred to as a reinforcement learning (RL) problem, where the goal is to learn an optimal policy from these observations without prior knowledge of the underlying dynamics.\footnote{In this article, unless explicitly emphasized as offline RL, RL refers to online RL.} Consequently, the agent needs to collect a large amount of samples by interacting with the environment to learn the best actions and maximize its long-term cumulative rewards. This characteristic of learning from samples is a fundamental facet of RL.
In~\appref{\ref{app:rl}}, we employ proximal policy optimization (PPO, \citealp{ppo}), a popular policy-based on-policy method, to assess the effectiveness of RL algorithms in solving optimal control problems. In line with our expectations, PPO's performance is substantially lower than direct policy optimization with fully-known dynamics due to its sample inefficiency. Additional details regarding the training process can be found in~\ref{app:exp_details}.

\section{Comparisons and A Unified Framework} \label{sec:method}

\subsection{Comparative Analysis} \label{sec:compare}

We summarize the characteristics of supervised learning and direct policy optimization approaches in Table~\ref{tab:sum}. In the training process, the supervised learning approach depends on an offline-generated dataset, while the batches in direct policy optimization are randomly sampled online. Two approaches are developed upon two different solution philosophies in optimal control: supervised learning is related to the open-loop optimal control, while direct policy optimization comes from the direct method for closed-loop optimal control.

\begin{table}[!htbp]
\small
\centering
\caption{Comparison on Supervised Learning (SL) and Direct Policy Optimization (DO)}
\begin{tabular}{ccccc}
\hline
\textbf{Methods}             & \textbf{Data} & \textbf{Training} & \textbf{Methodological Relevance}    &    \textbf{Challenge}               \\ \hline
\textbf{SL} &    \ding{51}         & Offline  & Available open-loop method  & \revisea{Dataset quality}  \\
\textbf{DO}                & \ding{56}             & Online            & Closed-loop direct method   & \revisea{Network optimization}   \\ \hline
\end{tabular}
\label{tab:sum}
\end{table}

In supervised learning, the information or the supervised signal comes from the open-loop optimal control solution, which can be obtained by any available numerical methods. 
Besides the aforementioned BVP or DDP solvers, we remark that the direct method can also be applied to open-loop problems, which is different from the direct policy optimization that solves the closed-loop OCP directly. It is often observed in practice that the most challenging part of supervised learning is to build an appropriate dataset using open-loop solvers rather than regression. It is why many efforts have been made, such as time-marching~\citep{Adaptive_BVP_HJB} or space-marching~\citep{ML_Enhanced_Landing} techniques, to improve the performance of open-loop solvers.
Furthermore, when the problem gets harder, the average solution time of each path gets longer. Given the same computation budget, one wishes to obtain a dataset of higher quality through adaptive sampling. 
One such improvement is related to data distribution, \textit{i.e.}, the discrepancy between the distribution of training data and the distribution of states indeed encountered by a controller over the feedback process increases over time, which is called the distribution mismatch phenomenon~\citep{long2022perturbational,IVP_Enhanced}. The IVP-enhanced sampling proposed by \cite{IVP_Enhanced} helps alleviate distribution mismatch and improves performance significantly. 
However, as shown in our numerical results, even with adaptive sampling
and sufficiently small validation error, 
there might still be a gap of realized total costs between the supervised learning outcome and optimal control in challenging problems.

Different from supervised learning, direct policy optimization needs no generated dataset yet has no prior information on optimal solutions. The whole problem is turned into a large-scale nonlinear optimization problem, which is natural to compute yet without further mathematical principles to follow. All difficulties are burdened by the optimization side, which makes it hard to train and leads to much longer training time. Consequently, direct policy optimization is much more time-consuming than the supervised learning approach. Relatedly, it is often observed that, given similar training time, supervised learning surpasses direct optimization in closed-loop simulation. What is worse, in challenging cases, the direct policy optimization may fail to get a reasonable solution under long-time training, since the randomly initiated policy network has a large deviation from the optimal solution and it is hard in such a complex optimization problem to find an appropriate direction and stepsize to converge to the optimal solution. 

\subsection{A Unified Training Paradigm}

\begin{figure}[!ht]
    \centering
    \includegraphics[width=0.973\textwidth]{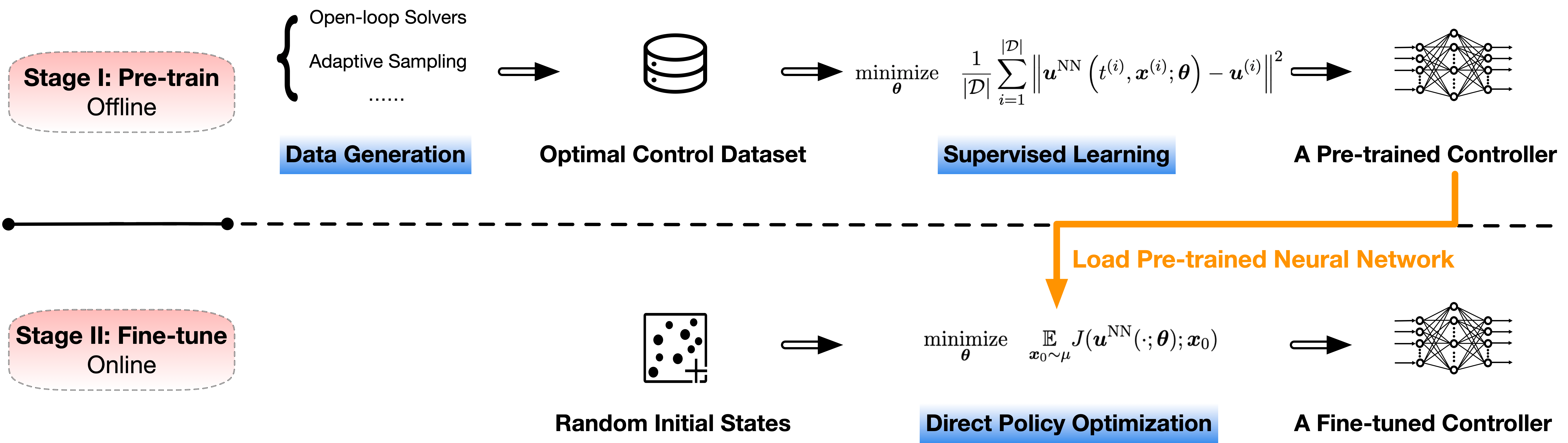}
     \caption{A unified training paradigm for neural network-based closed-loop optimal control consists of two stages. In Stage I, we first solve corresponding open-loop OCP to generate a dataset, on which we train the controller through supervised learning. In Stage II, we fine-tune the controller pre-trained in Stage I through online direct policy optimization.}
    \label{fig:pipeline}
\end{figure}

Summarizing the limitations discussed above, we find that supervised learning suffers from the dataset to further improve, while direct optimization needs proper initialization. Challenges in the two methods are different and somehow orthogonal, which offers us an opportunity to draw on their merits and complement each other. We combine the two methods into a unified training paradigm for neural network-based optimal control problems. The new paradigm, called \textbf{Pre-train and Fine-tune}, can be briefly sketched to pre-train by \textit{offline supervised learning} first and then fine-tune by \textit{online direct policy optimization}.
This training paradigm is outlined in Figure~\ref{fig:pipeline}, 
in which the two separated training approaches are combined as sequential training stages. In the first stage, we pre-train a controller via supervised learning. Various methods, such as adaptive sampling, can be included to achieve better performance. However, in some challenging cases, even with those techniques, there is still a gap between the learned controller and the optimal control. Thereby, we apply direct policy optimization based on the pertained network for fine-tuning. As discussed, a bottleneck of direct policy optimization is the extremely huge initial loss, which prevents the method to find a proper way to optimize. Fortunately, the pre-training through supervised learning can provide a reasonable solution that is close enough to the optimal solution. We can validate the optimality of different approaches in closed-loop simulation. Intuitively, the performances sorted from the worst to the best should be: \textit{Direct Policy Optimization $\leq$ Supervised Learning $\leq$ Pre-train and Fine-Tune $\leq$ Optimal Control}, which will be verified by experiments in Section~\ref{sec:exp}.

\revisea{The philosophy of \textit{Pre-train and Fine-tune} is relatively broad, encompassing various approaches, such as pre-training on a large dataset for a general task and fine-tuning on domain-specific data for more specialized tasks. In our context, it is important to emphasize that the objectives in the two stages are distinct yet correspond to the same optimal control problems from the same distribution of initial states. This approach bears some resemblance to certain combinations of offline and online RL algorithms~\citep{nair2020offline, lee2022offline, levine2020offline}. Offline RL involves learning from a fixed dataset, whereas online RL optimizes through dynamic interactions. A notable distinction in our method is its ability to generate a significant amount of optimal data when learning from a dataset of expert demonstrations. This is in contrast to many offline RL scenarios, where acquiring expert policy labels, such as human demonstrations, poses a challenge. Consequently, these scenarios often rely on static datasets, characterized by limited data and policies that may deviate from optimality.}
We believe that our Pre-train and Fine-tune strategy for optimal control problems offers a more systematic, controllable, and principled testbed for studying related algorithms, providing valuable insights into RL problems.

\section{Experiments} \label{sec:exp}

\subsection{Experimental Settings}

We consider two optimal control problems: the optimal attitude control of a satellite and the optimal landing of a quadrotor. \revisea{Both problems aim at controlling an aero-spatial vehicle from a general starting state to a target state with minimal cost, which has important engineering applications. Both problems have served as benchmarks in previous studies~\citep{ Adaptive_BVP_HJB, Satellite_17, IVP_Enhanced, Satellite_15, Quadrotor_04, Quadrotor_06, Quadrotor_12}. 
The optimal landing problem for a quadrotor is more challenging than that of the satellite due to the higher dimension and non-linear dynamics (see the detailed formulations of the problems in~\ref{app:dyn}). Besides, in~\ref{app:zero_control_lqr}, we evaluate the total cost of an uncontrolled system and the Linear Quadratic Regulator (LQR) controller as baselines to further underscore that these simpler controllers may have a big performance gap compared to the optimal controller.}

We focus on the problems where the initial states fall in a compact set $\mathcal{S}$ (defined in each example below respectively). In direct policy optimization, the initial state distribution $\mu$ in \eqref{eq3:obj-direct} is taken as the uniform distribution on $\mathcal{S}$. The dataset in supervised learning is also generated by uniformly sampling the initial states from $\mathcal{S}$ unless otherwise specified as the adaptive dataset in Section~\ref{sec:exp_quadrotor}.
For comprehensive comparisons, we conduct closed-loop simulations starting from different initial states. 
In robustness experiments, we corrupt the input states with uniform noises, \textit{i.e.}, the input for controllers becomes $\boldsymbol{x}(t) + \boldsymbol{n}(t)$, where each dimension of $\boldsymbol{n}(t)$ is independently and uniformly sampled from $[-\sigma, \sigma]$ with $\sigma$ being the disturbance scale.
As for evaluation metrics, we compute the pathwise ratio between the cost under the network controller and the optimal cost starting from the same initial state, called the \textbf{cost ratio}. \revisec{We remark that neither the BVP solution nor direct policy optimization can ensure global optimality. We use the cost ratio as a reasonable metric to establish a standard for comparing against a potentially optimal cost obtained by the open-loop solver.}
We plot the cumulative distribution curve and summarize different statistics of the cost ratio (with respect to the initial state distribution $\mu$) to evaluate different controllers.
All training details and time consumption are provided in \appref{\ref{app:exp_details}}. 

\subsection{The Optimal Altitude Control Problem of Satellite} \label{sec: exp_sat}

\paragraph{Settings} In this section, we conduct experiments on a satellite altitude control problem with a six-dimensional state~\citep{ Adaptive_BVP_HJB, Satellite_15, Satellite_17}. The state can be formulated as $\boldsymbol{x}=(\boldsymbol{v}^\mathrm{T}, \boldsymbol{\omega}^\mathrm{T})^\mathrm{T} = (\phi, \theta, \psi, \omega_1, \omega_2, \omega_3)^\mathrm{T}$, where $\boldsymbol{v}$ represents the attitude of the satellite and $\boldsymbol{\omega}$ represents the angular velocity. The full dynamics can be found in \appref{\ref{app:dyn_sat}}. The problem is to apply a torque $\boldsymbol{u} \in \mathbb{R}^3$ to stabilize the satellite to a final state of $\boldsymbol{v=0}$ and $\boldsymbol{\omega=0}$ at a fixed terminal time $T=20$. The set of interest for the initial state is $
\mathcal{S}_{\text{sate}}=\{\boldsymbol{v}, \boldsymbol{\omega} \in \mathbb{R}^3 \mid-\frac{\pi}{3} \leq \phi, \theta, \psi \leq \frac{\pi}{3}, -\frac{\pi}{4} \leq \omega_1, \omega_2, \omega_3 \leq \frac{\pi}{4}\}.
$

With the time-marching technique, it takes about 0.5 seconds on average to solve an open-loop solution for one trajectory and less than 1 minute in total (without parallel computing) to generate the whole optimal control dataset. Then we train the network under supervised learning for 100 epochs in 1 minute. In direct optimization, it costs more than 1 hour to train 2000 iterations from scratch, which is much longer than supervised learning.

\begin{table}[!ht]
\centering
\caption{Statistics of cost ratio and total computation time (min) in satellite's optimal attitude control problem}
\label{tab:deter_sat}
\begin{tabular}{cccccccc}
\hline
\textbf{Method}              & \textbf{Mean} & \textbf{Std} & \textbf{Max} & \textbf{Min} & \textbf{Median} & \textbf{Time}\\ \hline
\textbf{Direct Optimization}          &  1.048             &  0.034            &  1.180            &  1.012                             &   1.037   & 87           \\
\textbf{Supervised Learning}         & 1.003         & 0.001       & 1.006      & 1.002         & 1.003   & 2         \\
\hline
\end{tabular}
\end{table}

\begin{figure}[!ht]
    \centering
    \includegraphics[width=1.0\textwidth]{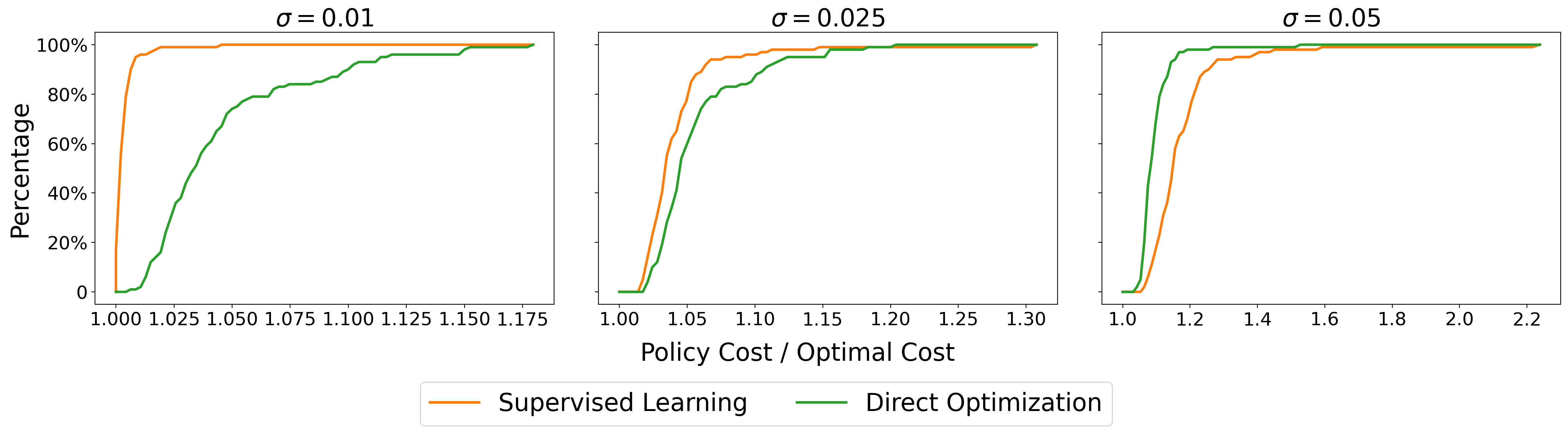}
    \caption{Cumulative distribution function of the cost ratio in the satellite problem with uniform disturbances of $\sigma=0.01, 0.025, 0.05$. 
    The spans of the horizontal axis are different. }
    \label{fig:sat_cum}
\end{figure}

\noindent \textit{Is supervised learning superior in performance?} We first evaluate the learned controllers in deterministic environments and summarize statistics of the cost ratio and total time in Table~\ref{tab:deter_sat}. As shown in Table~\ref{tab:deter_sat}, the performance of supervised learning is better than online direct policy optimization. The mean of cost ratio of supervised learning is extremely close to one, implying that the learned controller is close enough to the optimal control. Under such circumstances, there is no need for further fine-tuning. However, in the next challenging example, when purely supervised learning can not achieve competitive performance to optimal control, fine-tuning matters. \revisec{We further conduct experiments under different scales of disturbances $\sigma=0.01, 0.025, 0.05$ and plot the cumulative distribution functions of the cost ratio in Figure~\ref{fig:sat_cum}. 
The supervised controller can still stabilize the system and surpass direct optimization under small disturbances, demonstrating its robustness. We remark that since supervised learning is trained to fit a deterministic dataset rather than designed for stochastic systems, there may be states far away from those in the used optimal control dataset when the noise is large. Thus, supervised learning only maintains its priority to direct optimization under moderate disturbances.}

\begin{figure}[!ht]
    \centering
    \includegraphics[width=1.0\textwidth]{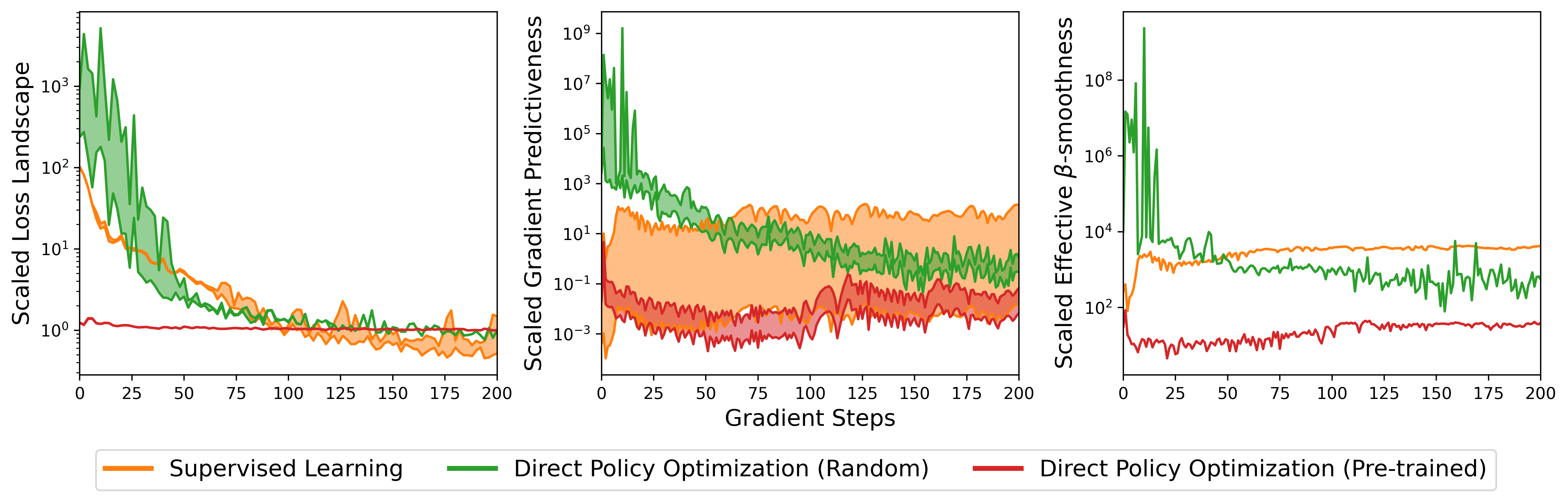}
  \caption{Comparative analysis of the optimization landscape in direct policy optimization and supervised learning. The figure displays, from left to right, normalized variations in different losses, normalized variations of $l_2$ changes in the gradient, and the normalized effective $\beta$-smoothness (the maximum ratio between gradient difference (in $l_2$-norm) and parameter difference) as moving in the gradient direction with different $step\_size$.}
    \label{fig:opt_landscape}
\end{figure}

\paragraph{What are the challenges in direct policy optimization? Exploring the optimization landscapes} In Figure~\ref{fig:opt_landscape}, we investigate the challenges associated with direct policy optimization by examining the optimization landscapes using three metrics inspired by~\cite{bnopt}. Let $l(\theta)$ denote the loss function (employed in supervised learning or direct policy optimization) and $lr$ represent the corresponding learning rate. To understand the landscape at various optimization stages given the current parameter $\hat{\theta}$, we assess local sensitivity by computing $l(\theta')$ and $\nabla_{\theta} l(\theta')$, where multiple $\theta'$ are given by $\theta^\prime = \hat{\theta} - {step\_size} \times \nabla_{\theta} l(\hat{\theta})$ with multiple ${step\_size}$, mirroring updates in stochastic gradient descent (SGD). Note that such local evaluations do not impact the training procedure itself and it is performed independently to gain insights into the behavior and characteristics of the optimization process. During local evaluations, the actual ${step\_size}$ is determined by scaling the learning rate used in training, denoted as $step\_size = step\_ratio \times lr$. This allows us to analyze the sensitivity to different learning rates along the gradient direction. In our implementation, we adopt multiple values of the $step\_ratio$ parameter within a predefined range and observe the corresponding changes in loss, which collectively form the loss landscape. Additionally, we compute the $l_2$ changes in the gradient to measure approximately the ``gradient predictiveness". We also compute the effective ``$\beta$-smoothness", as defined in \cite{bnopt}, which quantifies the maximum ratio between difference (in $l_2$-norm) in the gradient and the distance moved in a specific gradient direction as the $step\_ratio$ varies within the range. This measure provides valuable insights into the Lipschitz continuity of the gradient. \reviseb{To enable a consistent comparison between methods with different losses, we have normalized each one by scaling it according to the loss value in the final iteration, which brings metrics to a unified scale.
}

\reviseb{For both supervised learning and direct policy optimization (randomly initialized or pre-trained by supervised learning), Figure~\ref{fig:opt_landscape} displays scaled variations in losses (left), scaled variations of $l_2$ changes in the gradient (middle), and the scaled effective $\beta$-smoothness, the maximum ratio between gradient difference (in $l_2$-norm) and parameter difference (right), as moving in the gradient direction with different $step\_ratio$.\footnote{We remark that the valid domain of the $step\_ratio$ parameter used in direct policy optimization is considerably smaller compared to supervised learning since the former has a much rougher optimization landscape, as shown in Figure~\ref{fig:opt_landscape}. Further details regarding this observation can be found in \appref{\ref{app:exp_details}}.} 
Our comparisons shed light on the challenges associated with direct policy optimization, particularly its sensitivity to the actual $step\_size$.
The much higher shaded regions in Figure~\ref{fig:opt_landscape} (left and middle) for the direct policy optimization (randomly initialized) show its larger sensitivity to the $step\_size$ compared to supervised learning. 
Furthermore, the effective $\beta$-smoothness, which measures the Lipschitzness of the gradient, demonstrates that the optimization landscape in supervised learning is notably smoother. Significantly, the metrics for direct policy optimization which is pre-trained by supervised learning, exhibit the smoothest outcomes in the comparison, showing the critical role of proper initialization. These findings emphasize the challenges faced in direct policy optimization and highlight the advantages of supervised learning, which benefits from a smoother and more benign optimization landscape.}

\subsection{The Optimal Landing Problem of Quadrotor} \label{sec:exp_quadrotor}
\paragraph{Settings} We consider a more complex twelve-dimensional problem aiming to control a quadrotor to land at a target position in a fixed terminal time $T$~\citep{IVP_Enhanced}. The state can be formulated as $\boldsymbol{x} = \left(\boldsymbol{p}^{\mathrm{T}}, \boldsymbol{v}_b^{\mathrm{T}}, \boldsymbol{\eta}^{\mathrm{T}}, \boldsymbol{w}_b^{\mathrm{T}}\right)^{\mathrm{T}}$, where $\boldsymbol{p}=(x, y, z)^{\mathrm{T}} \in \mathbb{R}^3$ denotes the position in Earth-fixed coordinates, $\boldsymbol{v}_b=(v_x, v_y, v_z)^{\mathrm{T}} \in \mathbb{R}^3$ denotes the velocity with respect to the body frame, $\boldsymbol{\eta}=(\phi, \theta, \psi)^{\mathrm{T}} \in \mathbb{R}^3$ denotes the attitude in Earth-fixed coordinates, and $\boldsymbol{w}_b \in \mathbb{R}^3$ denotes the angular velocity in the body frame. The full dynamics can be found in \appref{\ref{app:dyn_qua}}. We train a feedback controller $\boldsymbol{u}^\text{NN}(t, \boldsymbol{x}): [0, T] \times \mathbb{R}^{12} \rightarrow \mathbb{R}^4$ to control the rotor thrusts to land the quadrotor.
The detailed setting is similar to that in~\cite{IVP_Enhanced}. Specifically, the set of interest for the initial state is $\mathcal{S}_{\text{quad}} = \{x, y \in [-40, 40], z \in [20, 40],~~v_x, v_y, v_z \in [-1, 1], \theta, \phi \in [-\pi/4, \pi/4], \psi \in [-\pi, \pi], \boldsymbol{w} = \boldsymbol{0} \}$.
We also consider problems whose initial states lie in a smaller domain $\tilde{\mathcal{S}}_{\text{quad}} \subset \mathcal{S}_{\text{quad}}$ to study the effect of varying time horizons, where $\tilde{\mathcal{S}}_{\text{quad}} = \{x, y \in [-8, 8], z \in [4, 8], v_x, v_y, v_z \in [-0.2, 0.2], \theta, \phi \in [-\pi/20, \pi/20], \psi \in [-\pi/5, \pi/5], \boldsymbol{w} = \boldsymbol{0} \}$. %

\begin{figure}[!ht]
    \centering
    \includegraphics[width=1.0\textwidth]{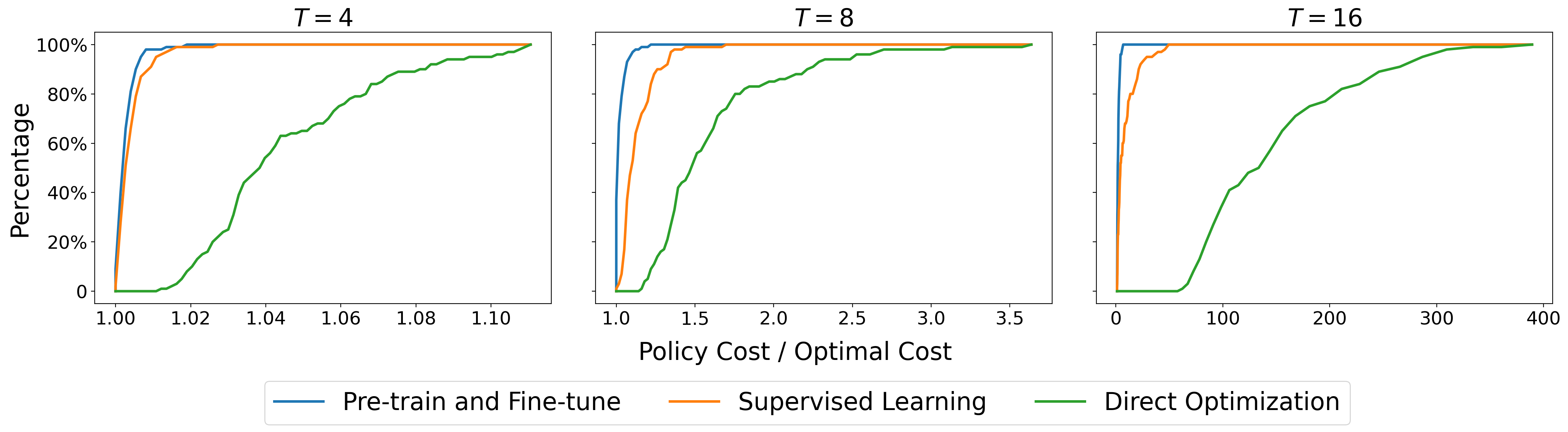}
  \caption{Cumulative distribution function of the cost ratio in quadrotor's optimal landing problem on $x_0 \in \tilde{\mathcal{S}}_{\text{quad}}$ with varying time horizons $T=4, 8, 16$. Note that the spans of the horizontal axis increase from left to right.}
    \label{fig:qua_medium_time}
\end{figure}

\paragraph{What are the challenges in direct policy optimization? Comparisons on different horizons} As discussed in Section~\ref{sec:method}, the challenge in direct policy optimization is the optimization itself. To better understand this challenge, we first consider the initial state $\boldsymbol{x}_0$ uniformly sampled from $\tilde{\mathcal{S}}_{\text{quad}}$ and the time horizon $T$ equals 4, 8, or 16. The cumulative distribution functions of cost ratio are plotted in Figure~\ref{fig:qua_medium_time} and the corresponding statistics are reported in Table~\ref{tab:varying_time}. 
From the results of the direct optimization, we can clearly observe that, as the total time $T$ gets longer, the method deteriorates dramatically, and the gap between direct optimization and supervised learning grows significantly. The reason is when the horizon is longer, the forward trajectory controlled by a randomly initialized network will severely deviate from the optimal path more and more, resulting in a larger initial loss which grows quickly as $T$ increases. Under such circumstances, the optimization is much harder since finding a proper direction for the neural network to improve is difficult. Compared to training from scratch, the significant improvements brought by a pre-trained network emphasize the main challenge in direct optimization is the optimization itself. % 

\begin{table}[!ht]
\caption{Statistics of cost ratio and total computation time (min) in quadrotor's optimal landing problem on $\boldsymbol{x}_0 \in \tilde{\mathcal{S}}_{\text{quad}}$ with different time horizons $T=4, 8, 16$.} 
\label{tab:varying_time}
\centering
\begin{tabular}{cccccccc}
\hline
\multicolumn{1}{l}{\textbf{Horizon}} & \textbf{Method}               & \textbf{Mean}        & \textbf{Std}         & \textbf{Max}         & \textbf{Min}     & \textbf{Median}   & \textbf{Time}      \\ \hline
\multirow{3}{*}{\textbf{$\bm{T}$=4}}        & \textbf{DO}   & 1.05                   & 0.02                & 1.11                & 1.01                          & 1.04     & 106           \\
                                     & \textbf{SL}   & 1.01               & 0.00               & 1.03              & 1.00                   & 1.00    & 16            \\
                                     & \textbf{Fine-tune} & 1.00                & 0.00                 & 1.02                 & 1.00                            & 1.00       & 16 + 4          \\ \hline
\multirow{3}{*}{\textbf{$\bm{T}$=8}}        & \textbf{DO}   & 1.63                    & 0.43                 & 3.70                 & 1.17                     & 1.50    & 157                \\
                                     & \textbf{SL}   & 1.15                 & 0.11                 & 1.72                 & 1.01                              & 1.11   & 16           \\
                                     & \textbf{Fine-tune} & 1.03                 & 0.03                 & 1.23                 & 1.00                            & 1.02   & 16 + 6              \\ \hline
\multirow{3}{*}{\textbf{$\bm{T}$=16}}       & \textbf{DO}   & 157.71                    & 75.65                 & 420.18                 & 64.64                 & 143.10       & 387          \\
 & \textbf{SL}   & 9.59               & 11.00               & 51.18             & 1.11                          & 4.52    & 16          \\
                                     & \textbf{Fine-tune} & 2.31               & 1.22                 & 6.93                & 1.11                             & 1.83     & 16 + 7            \\ \hline
\end{tabular}
\end{table}

\begin{figure}[!ht]
    \centering
    \includegraphics[width=1.0\textwidth]{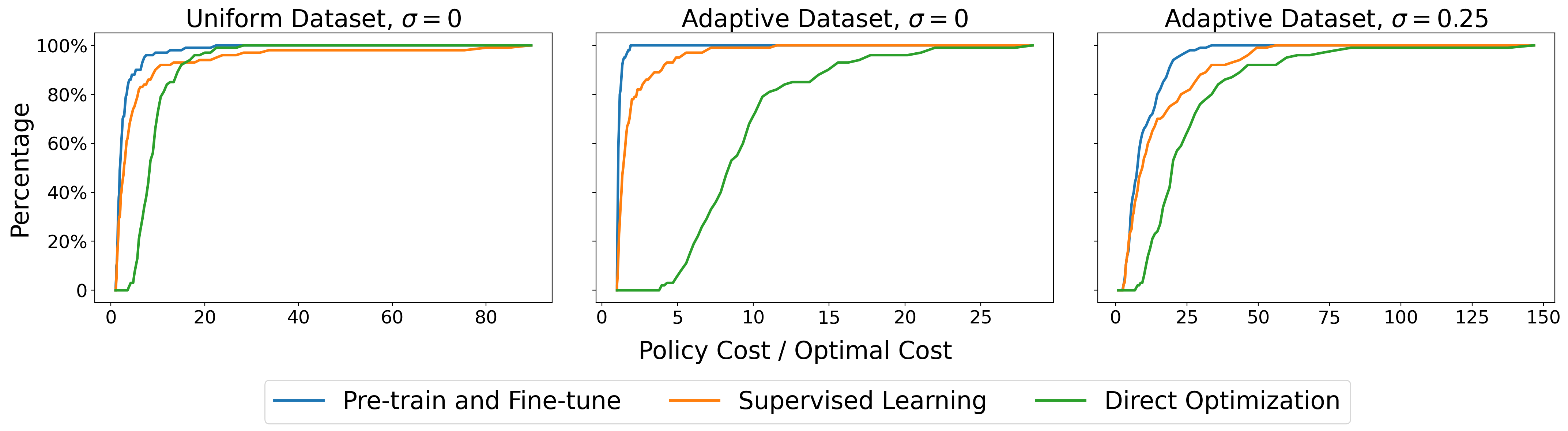}
  \caption{Cumulative distribution function of the cost ratio in quadrotor's optimal landing problem on $\bm{x}_0 \in \mathcal{S}_{\text{quad}}$ with $T=16$. Left and middle: closed-loop evaluation in deterministic environments where supervised learning is trained on the dataset with uniformly sampled initial states and the adaptive dataset generated by IVP-enhanced sampling respectively. Right: closed-loop evaluation in stochastic environments with $\sigma=0.25$ where supervised learning is trained on the adaptive dataset. The spans of the horizontal axis are different.}
  \label{fig:10_16_dis}
\end{figure}

\noindent \textit{What are the challenges in supervised learning?} Table~\ref{tab:exp_10_16} reports results in a more difficult setting where the initial state $\boldsymbol{x}_0$ is uniformly sampled from a larger set $\mathcal{S}_{\text{quad}}$ and the total time horizon $T=16$. Due to the problem's difficulty, we consider two types of datasets of the same size in supervised learning, one from uniform sampling and one from IVP-enhanced adaptive sampling~\citep{IVP_Enhanced}, and we apply fine-tuning in both cases. Supervised learning surpasses the direct method still, yet it has poor performance in the worst cases, for example, the max cost ratio is very large. This is due to the distribution mismatch phenomenon~\citep{IVP_Enhanced}: when the controller encounters states far from those in the optimal control dataset, the control output is unreliable.\footnote{The exceptionally large cost ratio in the worst cases, due to the distribution mismatch phenomenon, also causes the mean of the cost ratio to fluctuate significantly across multiple runs of experiments when learning from a fixed dataset. In Table \ref{tab:exp_10_16}, we report the best mean among multiple runs. Meanwhile, the median of the cost ratios across different runs fluctuates much less and consistently reflects the superiority of supervised learning over the direct method.}
To alleviate this phenomenon, we employ IVP-enhanced sampling to obtain an adaptive dataset for supervised learning.
We remark that the validation loss of supervised learning has similar behavior and converges to a small scale of \texttt{1e-4} on both the uniform dataset and the adaptive dataset, but the closed-loop simulation performances are much different as shown in Table~\ref{tab:exp_10_16} and Figure~\ref{fig:10_16_dis}. This fact supports our analysis in Section~\ref{sec:compare} on supervised learning whose challenge is the dataset instead of the optimization process. Furthermore, comparisons between supervised learning and direct policy optimization in Figures~\ref{fig:qua_medium_time}, \ref{fig:10_16_dis} and Tables~\ref{tab:varying_time}, \ref{tab:exp_10_16} provide similar evidence as in the example of the satellite: supervised learning achieves better results than direct policy optimization with less time, which also implies that the learning procedure in supervised learning is much easier than direct policy optimization. %

\begin{table}[!ht]
\caption{Statistics of cost ratio and total computation time (min) in quadrotor's optimal landing problem on $\boldsymbol{x}_0 \in \mathcal{S}_{\text{quad}}$ and $T=16$. The pre-training is conducted with a uniformly sampled dataset and an adaptively sampled dataset, respectively.}
\label{tab:exp_10_16}
\centering
\begin{tabular}{clccccccc}
\hline
\textbf{Method}            & \textbf{Fine-tune} & \textbf{Mean}        & \textbf{Std}         & \textbf{Max}         & \textbf{Min}        & \textbf{Median}      & \textbf{Time}   \\ \hline
\textbf{Direct}                        & -----                  & 9.72                    & 4.31                & 29.67                & 3.96                  & 8.69    & 426            \\ \hline
\multirow{3}{*}{{\begin{tabular}{c}
                   \textbf{Pre-train}\\
                    (Uniform)
                   \end{tabular}}}      & 0 (SL)                  & 6.45                    & 13.18               & 94.96              & 1.11                     & 2.97 & 68                 \\
                                     & 100                & 3.07                    & 3.11                 & 23.22                 & 1.06                      & 2.09  & 68 + 9                \\
                                     & 1000               & 1.55                    & 0.69                 & 6.89                 & 1.05                           & 1.34  & 68 + 92                 \\ \hline
\multirow{3}{*}{{\begin{tabular}{c}
                   \textbf{Pre-train}\\
                    {(Adaptive)}
                   \end{tabular}}}       & 0 (SL)                & 2.05                    & 1.61                 & 11.84                 & 1.03           & 1.48      & 225           \\
                                     & 100                & 1.26                    & 0.34                 & 3.10                 & 1.03                 & 1.13   &   225 + 9            \\
                                     & 1000               & 1.06                    & 0.04                 & 1.28                 & 1.01                  & 1.05   & 225 + 92              \\ \hline
\end{tabular}
\end{table}

\noindent \textit{Can Pre-train and Fine-tune improve performance and robustness?} In all settings reported in Table~\ref{tab:varying_time} and Table~\ref{tab:exp_10_16}, fine-tuning improves the performance of supervised learning as desired, regardless of whether the dataset in supervised sampling is uniformly or adaptively sampled. Figure~\ref{fig:10_16_dis} plots the cumulative distribution functions of cost ratio when $\boldsymbol{x}_0 \in \mathcal{S}_{\text{quad}}$ and the measurement noise $\sigma=0$ or $0.25$.
More experiments under different measurement noises are provided in \appref{\ref{app:results_sto}}, demonstrating similar improvements brought by fine-tuning on robustness. We highlight the fine-tuning time (see Table~\ref{tab:varying_time}, ~\ref{tab:exp_10_16} and \appref{\ref{app:training_time}}), which is only about a few minutes with significant improvements in both performance and robustness. Direct policy optimization costs several hours and achieves sub-optimal solutions when training from scratch, whilst fine-tuning can improve performance significantly within a few minutes based on a proper initialization. 
{We stop the pre-training process upon observing little progress in the validation loss, followed by initiating the fine-tuning phase.} As shown in Table~\ref{tab:exp_10_16},
it first takes about 1 hour to pre-train and then 10 minutes to improve the performance significantly through fine-tuning 100 iterations, while the direct optimization requires about 7 hours to train from scratch. If we further fine-tune the controller to 1000 iterations, it will cost more than 1 hour yet improve less than before. Therefore, the results of the Pre-train and Fine-tune strategy we report are fine-tuned by 100 iterations, unless otherwise stated in Table~\ref{tab:exp_10_16}.

\begin{figure}[!ht]
    \centering
    \includegraphics[width=1.0\textwidth]{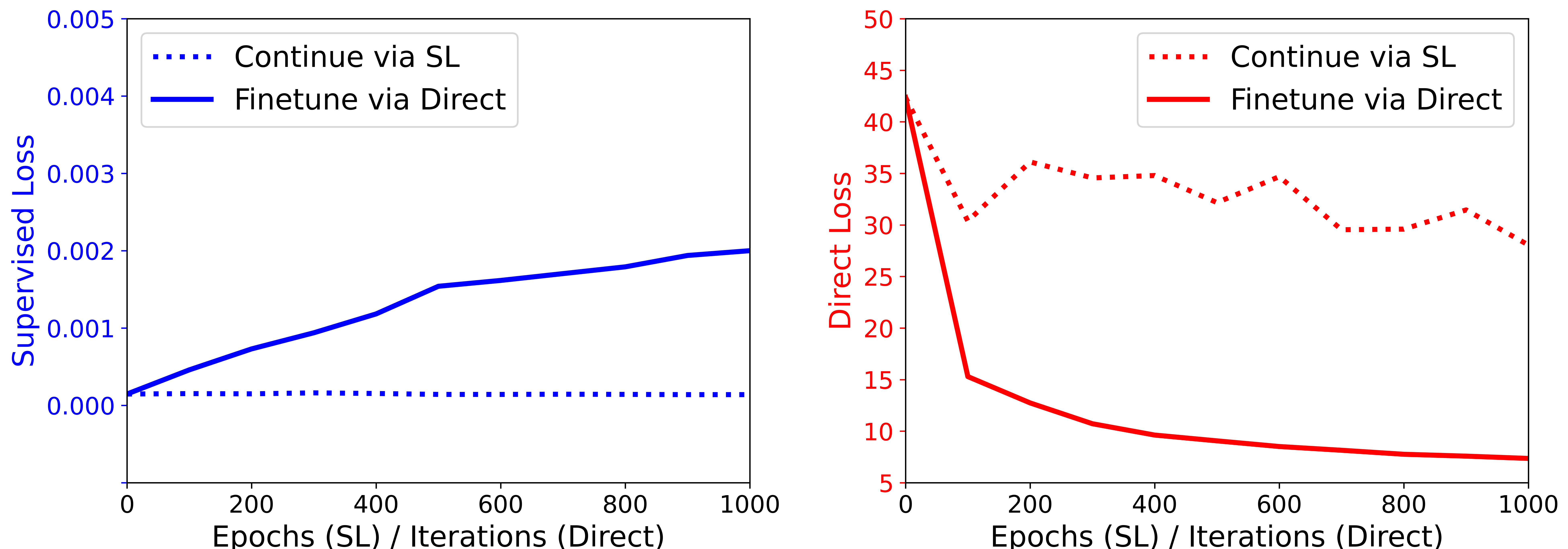}
  \caption{Loss curves for continuing via supervised learning and fine-tuning via direct policy optimization based on the same pre-trained model.}
  \label{fig:further_train}
\end{figure}

We proceed with additional experiments to highlight the distinction between further training via supervised learning and fine-tuning via direct optimization. Using the same pre-trained model from the second line Table~\ref{tab:exp_10_16}, we conduct additional SL training for 1000 epochs and fine-tuning for 1000 iterations, respectively. Figure~\ref{fig:further_train} visually depicts the different loss behaviors observed during further training with different loss functions. Continuing with SL yields only marginal improvements in the direct loss (average of the total costs for the samples). Conversely, fine-tuning through direct optimization results in rapid improvements in the direct loss, while even having a litter adverse impact on the SL loss. %  
This comparison shows that there can be two models which perform similarly on the SL loss but significantly different on the optimal control objective. In other words, for models that are reasonably close to the optimal control, the SL loss may not be an effective indicator of how the model solves the optimal control problem.
Therefore, we emphasize the significance of fine-tuning when the improvements achieved through SL training are marginal and fail to substantially enhance the overall performance.

\section{Conclusion}

In this work, we conduct a comprehensive comparative study of two approaches for training a neural network-based closed-loop optimal controller. Primarily, we establish a benchmark for comparing offline supervised learning and online direct policy optimization and analyze the merits and drawbacks of the two methods in detail. The experimental results highlight the priority of offline supervised learning on both performance and training time. We point out that the main challenges to the two methods are dataset and optimization, respectively. Based on the detailed analysis, we naturally propose the Pre-train and Fine-tune strategy as a unified training paradigm for closed-loop control, which significantly enhances performance and robustness.

\appendix

\section{Additional Results}

\subsection{Direct Policy Optimization with Fully-Known Dynamics} \label{app:adjoint}

As mentioned in Section~\ref{sec: direct}, the objective (\ref{eq3:obj-direct}) of direct policy optimization with fully-known dynamics can be optimized through stochastic gradient descent. There are two ways to compute the required (stochastic) gradient: back-propagation and the adjoint method.

In back-propagation~\citep{bp}, we first discretize the continuous trajectory with a given initial state $\boldsymbol{x}_0$ by an ODE integration scheme and get a discrete trajectory to approximate the total cost $J(\boldsymbol{u}; \boldsymbol{x}_0)$. 
Based on discrete trajectories, we apply back-propagation through the operations of the trajectories to compute the gradient with respect to the parameters $\boldsymbol{\theta}$. 
Adaptive ODE solvers, which automatically find suitable stepsize under certain error tolerance, can often achieve a better approximation efficiency of the continuous trajectory than a fixed-stepsize integrator. However, directly leveraging adaptive ODE solvers leads to an implicit and unknown depth of the computation graph, resulting in a deterioration of performance and computation time in back-propagation. To overcome this limitation, \cite{ACA_NODE} suggest deleting the computation graph of adaptive solvers and storing the final computed stepsize only, which reduces the memory cost and benefits in the final performance.

The adjoint method~\citep{pontryaginMathematicalTheoryOptimal1962} takes another avenue to calculate the gradient of the total loss $J$ with respect to parameters $\boldsymbol{\theta}$ in a continuous viewpoint. \cite{CTPG} defines the adjoint $\boldsymbol{\alpha}(t)=\partial J / \partial \boldsymbol{x}(t)$ with an augmented ODE $\frac{d \boldsymbol{\alpha}(t)}{d t}=-\boldsymbol{\alpha}(t)^{\top} \frac{\partial \boldsymbol{f}(\boldsymbol{x}, \boldsymbol{u})}{\partial \boldsymbol{x}} - \frac{\partial L(\boldsymbol{x}, \boldsymbol{u})}{\partial \boldsymbol{x}}$, where the notations are the same as OCP (\ref{eq1-objective}). Then the gradient can be calculated by $\frac{\partial J}{\partial \boldsymbol{\theta}}=\int_{0}^{T} \boldsymbol{\alpha}(t)^{\top} \frac{\partial \boldsymbol{f}}{\partial \boldsymbol{u}}\frac{\partial \boldsymbol{u}}{\partial \boldsymbol{\theta}} + \frac{\partial L}{\partial \boldsymbol{u}}\frac{\partial \boldsymbol{u}}{\partial \boldsymbol{\theta}} d t$, which requires a forward ODE to solve the state $\boldsymbol{x}(t), t\in[0,T]$ and a backward ODE with $\boldsymbol{\alpha}(T)=\partial M / \partial \boldsymbol{x}(T)$ to solve the adjoint $\boldsymbol{\alpha}(t)$.
The numerical error and memory usage of the adjoint method can be less than the back-propagation; however, the computing time is longer due to the additional calls on ODE solvers. The idea of the adjoint method has been followed by Neural ODE~\citep{NeuralODE} to view the discrete neural network model as a continuous flow. Extending the implementation of~\cite{NeuralODE} back to control problems, efforts are made in neural network-based control under deterministic dynamics~\citep{AIPontryagin, CTPG}. Unfortunately, naively leveraging the adjoint method may cause catastrophic divergence, even in the simple linear quadratic regulator (LQR) case~\citep{CTPG}, which can be alleviated by adding checkpoints to ODE trajectories~\citep{CTPG, ACA_NODE, ANODE}. The aforementioned idea of adaptive solvers can also be applied and the stored states can serve as the checkpoints~\citep{ACA_NODE}.

We compare direct policy optimization with back-propagation or the adjoint method in the optimal landing problem on $\tilde{\mathcal{S}}_{\text{quad}}$ with $T=4, 8, 16$ in Table~\ref{tab:adjoint}. The adjoint method costs time about twice as back-propagation. We round up to 4 decimal places to show the slight differences in $T=4, 8$. When the time horizon is $T=16$, the adjoint method deteriorates even with checkpoints added. Since memory is not the bottleneck, we use back-propagation with competitive performances and faster training speed to train the direct optimization in the main content.

% % Please add the following required packages to your document preamble:
% % \usepackage{multirow}
% \begin{table}[!htbp]
% \caption{Cost ratio in quadrotor’s optimal landing problem on $x_0 \in \tilde{\mathcal{S}}_{quad}$ with varying time horizons $T = 4, 8, 16$ and comparisons on back-propagation (bp) and the adjoint method.}
% \label{tab:adjoint}
% \begin{tabular}{ccccccc}
% \hline
% \textbf{Problem}               & \textbf{Method}  & \textbf{Training Time (min)} & \textbf{Mean} & \textbf{Std} & \textbf{Max} & \textbf{Min} \\ \hline
% \multirow{2}{*}{\textbf{T=4}}  & \textbf{BP}      & 87                    & 1.0468        & 0.0235       & 1.1120       & 1.0125       \\
%                                & \textbf{Adjoint} & 128                    & 1.0449        & 0.0232       & 1.1089       & 1.0123       \\ \hline
% \multirow{2}{*}{\textbf{T=8}}  & \textbf{BP}      & 152                    & 1.6295        & 0.4349       & 3.6979       & 1.1712       \\
%                                & \textbf{Adjoint} & 298                    & 1.6218        & 0.4416       & 3.6976       & 1.1609       \\ \hline
% \multirow{2}{*}{\textbf{T=16}} & \textbf{BP}      &     329                   & 157.7106      & 75.6489      & 420.1823     & 64.6441      \\
%                                & \textbf{Adjoint} & 706                    & 174.5295      & 116.2190     & 787.1036     & 61.1059      \\ \hline
% \end{tabular}
% \end{table}

\begin{table}[!htbp]
\caption{Cost ratio and training time (min) in quadrotor’s optimal landing problem on $\boldsymbol{x}_0 \in \tilde{\mathcal{S}}_{\text{quad}}$ with varying time horizons $T = 4, 8, 16$ and comparisons on back-propagation (BP) and the adjoint method.}
\label{tab:adjoint}
\centering
\begin{tabular}{ccccccc}
\hline
\textbf{Problem}               & \textbf{Method}  & \textbf{Time} & \textbf{Mean} & \textbf{Std} & \textbf{Max} & \textbf{Min} \\ \hline
\multirow{2}{*}{\textbf{$\bm{T}$=4}}  & \textbf{BP}      & 106                    & 1.0468        & 0.0235       & 1.1120       & 1.0125       \\
                               & \textbf{Adjoint} & 207                    & 1.0449        & 0.0232       & 1.1089       & 1.0123       \\ \hline
\multirow{2}{*}{\textbf{$\bm{T}$=8}}  & \textbf{BP}      & 157                    & 1.6295        & 0.4349       & 3.6979       & 1.1712       \\
                               & \textbf{Adjoint} & 298                    & 1.6218        & 0.4416       & 3.6976       & 1.1609       \\ \hline
\multirow{2}{*}{\textbf{$\bm{T}$=16}} & \textbf{BP}      &     387                   & 157.7106      & 75.6489      & 420.1823     & 64.6441      \\
                               & \textbf{Adjoint} & 706                    & 174.5295      & 116.2190     & 787.1036     & 61.1059      \\ \hline
\end{tabular}
\end{table}

\subsection{Direct Policy Optimization with Reinforcement Learning} \label{app:rl}
In this subsection, we compare the performance of direct policy optimization with and without fully-known dynamics.  We focus on the optimal attitude control problem of a satellite. With fully-known dynamics, we optimize by back-propagating through the computational graph of the entire trajectory. With unknown dynamics, we utilize the on-policy, model-free method, Proximal Policy Optimization (PPO)~\citep{ppo}, as a representative method in reinforcement learning, which has been empirically demonstrated to achieve state-of-the-art results in a wide range of tasks, including robotics~\citep{ppo_robot} and large language models~\citep{instrcutgpt}.
\begin{table}[!ht]
\centering
\caption{Statistics of cost ratio in satellite's optimal attitude control problem}
\label{tab:rl_sat}
\begin{tabular}{ccccccc}
\hline
\textbf{Method}              & \textbf{Mean} & \textbf{Std} & \textbf{Max} & \textbf{Min} & \textbf{Median} \\ \hline
\textbf{BP With Fully-Known Dynamics}          &  1.048             &  0.034            &  1.180            &  1.012                             &   1.037           \\
\textbf{RL (1x Samples)}         & 12.60         & 8.959       & 78.52      & 4.106         & 10.31          \\
\textbf{RL (5x Samples)}         & 1.628         & 0.392       & 4.832      & 1.254         & 1.548          \\
\hline
\end{tabular}
\end{table}

\begin{figure}[!ht]
    \centering
    \includegraphics[width=1.0\textwidth]{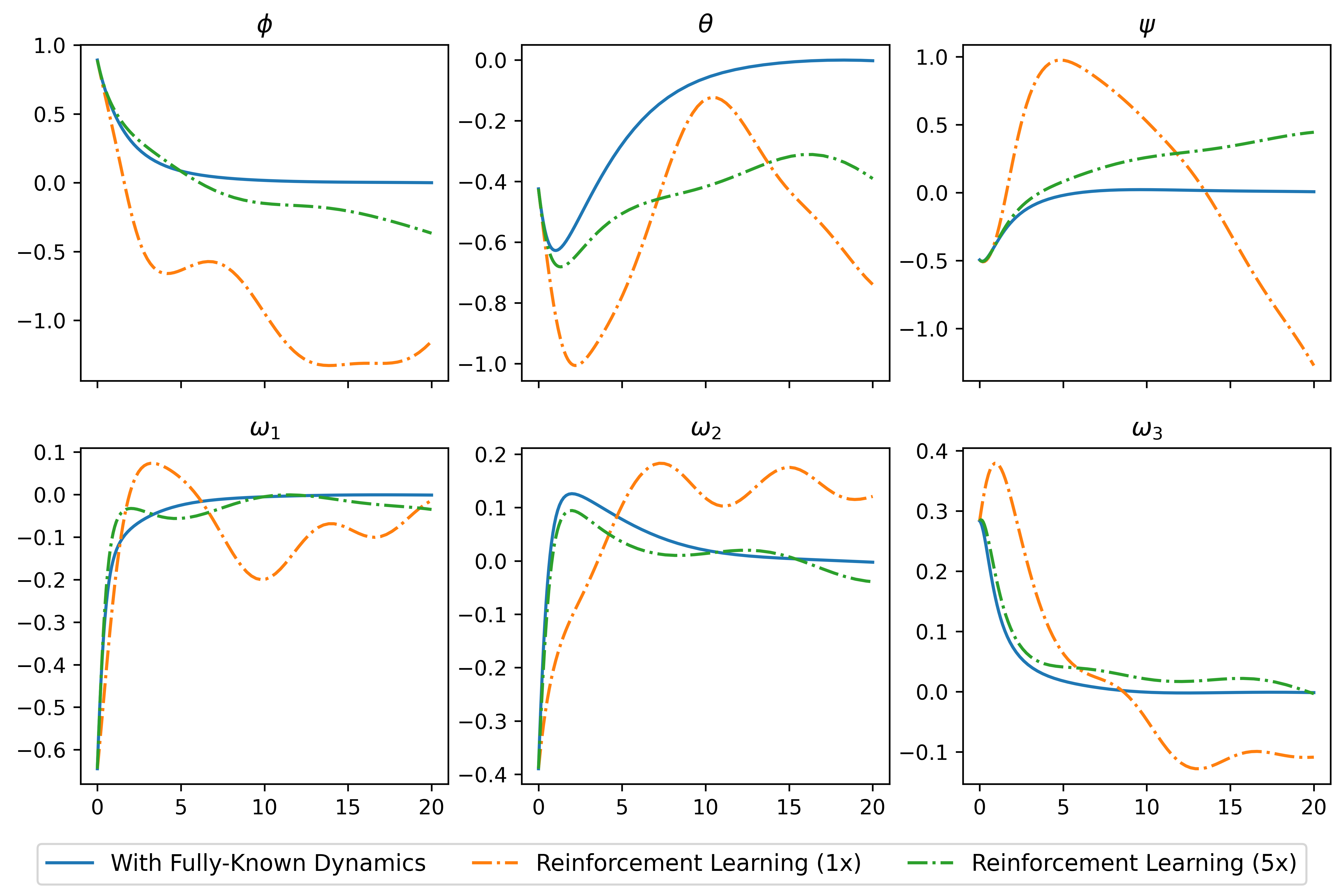}
  \caption{Each dimension of solutions by direct policy optimization with fully-known dynamics and reinforcement learning. Different sample scales used for reinforcement learning are denoted as 1x and 5x.}
    \label{exp:rl}
\end{figure}

To ensure a fair comparison, we conduct experiments with RL using either the same number of trajectories as direct policy optimization with fully-known dynamics or five times as many, denoted as \textbf{1x and 5x samples}, respectively. As presented in Table~\ref{tab:rl_sat}, direct policy optimization with fully-known dynamics has a cost ratio of $1.048 \pm 0.034$, which is comparable to optimal control. In contrast, RL (5x Samples) has a worse cost ratio of $1.423\pm0.434$, and RL (1x Samples) is significantly worse of $12.58 \pm 8.959$. These results illustrate that RL algorithms are significantly less sample-efficient.

We demonstrate the trajectories controlled by networks trained using the above methods in Figure~\ref{exp:rl}, starting from the same initial point. While the RL solutions seem plausible, they are less optimal than the ones learned through dynamics. In this case, the optimal cost is 4.03, and cost of direct optimization with fully-known dynamics is 4.06, which are very close. However, the cost of RL (5x) is 5.75 and the cost of RL (1x) is much worse, 27.20. The advantages of using fully-known dynamics are straightforward, as we possess more information about the trajectories. Further details of our experiments are provided in \ref{app:exp_details}.

% \yz{rl cost: 27.20
% rl bst cost: 5.75
% sl cost: 4.05
% BVP cost: 4.03, The BVP and direct almost the same so we only present the direct trajectory.}

\subsection{Robustness Evaluation} \label{app:results_sto}

%For robustness evaluation, 

\revisec{In this section, we present robustness experiments in quadrotor's optimal landing problem, when there is state measurement noise in the close-loop simulation.} We only fine-tune 100 iterations within 10 minutes based on the network pre-trained by supervised learning.

\begin{figure}[!htbp]
    \centering
    \includegraphics[width=1.0\textwidth]{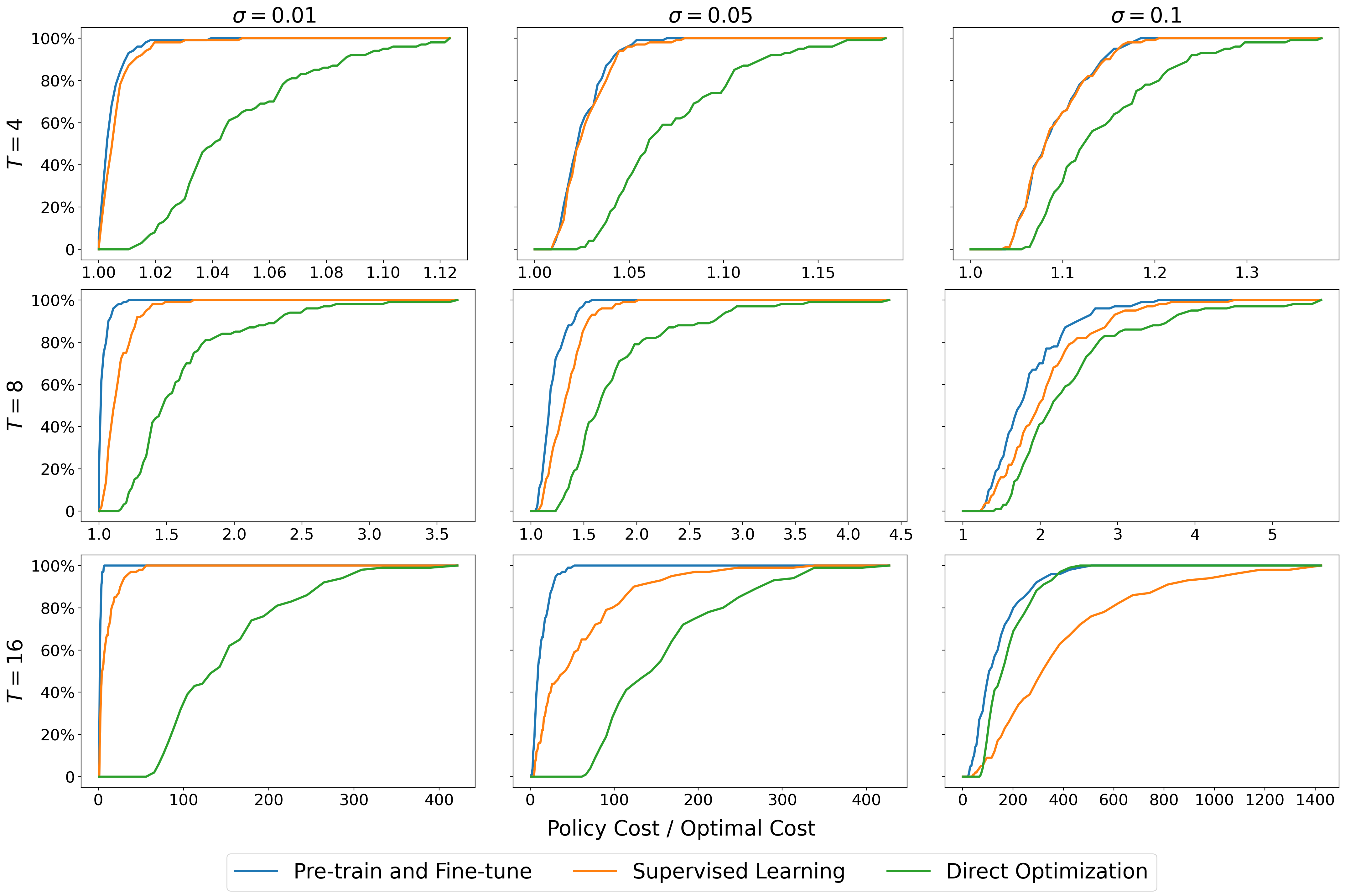}
  \caption{Cumulative distribution function of the cost ratio in quadrotor's optimal landing problem with $\boldsymbol{x}_0 \in \tilde{\mathcal{S}}_{\text{quad}}$. The rows from top to bottom are varying time horizons $T=4,8,16$. The columns from left to right are under different scales of disturbance $\sigma=0.01,0.05,0.1$.}
  \label{exp:robust_2_time_varying}
\end{figure}

\begin{figure}[!htbp]
    \centering
    \includegraphics[width=1.0\textwidth]{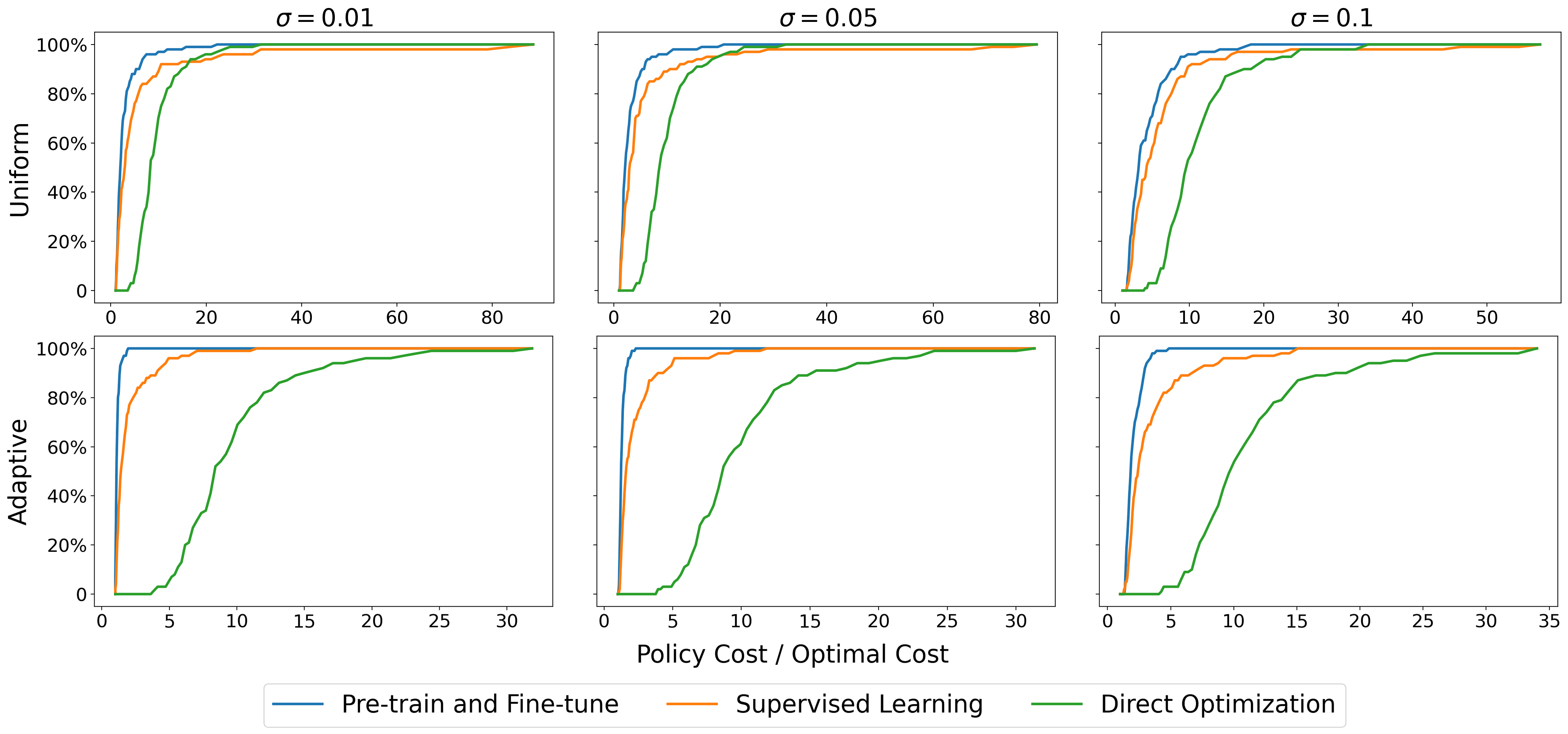}
  \caption{Cumulative distribution function of the cost ratio in quadrotor's optimal landing problem with $\boldsymbol{x}_0 \in \mathcal{S}_{\text{quad}}$ and $T=16$. The rows from top to bottom are related to different distributions of the dataset, uniform and adaptive. The columns from left to right are under different scales of disturbance $\sigma=0.01,0.05,0.1$.}
    \label{exp:robust_10_16}
\end{figure}

\revisec{Figure~\ref{exp:robust_2_time_varying} and \ref{exp:robust_10_16} demonstrate that supervised learning outperforms direct policy optimization under moderate disturbances. Additionally, we highlight the superior robustness of the Pre-train and Fine-tune strategy compared to both direct policy optimization and supervised learning.} When the noise is small, supervised learning is still able to control the quadrotor and maintain its priority over direct optimization. However, since supervised learning is not designed for stochastic optimal control, it does not suit the case when the noise is very large. Also, the long time horizon enlarges the influence brought by disturbance, accumulates the error along the trajectory, and results in the distribution mismatch phenomenon and poor performances. Fine-tuning brings robustness to the pre-trained network, which can endure much larger disturbance and performs the best even when the performance of supervised learning decreases dramatically as the bottom row in Figure~\ref{exp:robust_2_time_varying} shows. Figure~\ref{exp:robust_10_16} also shows that fine-tuning is able to improve the robustness regardless of the dataset distribution where it is pre-trained on.

\revisec{
\subsection{Uncontrolled System and LQR Performances}\label{app:zero_control_lqr}}
\revisec{
In this subsection, we aim to quantify the difficulty of different problems by evaluating the performance of various controllers. We begin by computing the total cost of uncontrolled systems on the test dataset. To provide a benchmark for comparison, we also consider the Linear Quadratic Regulator (LQR) as a baseline. The LQR controller is designed to optimize the system's performance by linearizing the full dynamics around the terminal state~\citep{mehrmann1991autonomous}. In the context of satellite altitude control, the terminal point is defined as $\boldsymbol{x_f}=0$ and $\boldsymbol{u_f}=0$, while in the quadrotor's landing problem, it is $\boldsymbol{x_f}=0$ and $\boldsymbol{u_f}=mg$. We apply Drake~\citep{drake} to solve the related finite horizon LQR problems. Following the same problem setup as discussed in the main content, the feedback control $\boldsymbol{u}(t, \boldsymbol{x})=-K(t)\boldsymbol{x}$ is utilized, where $K(t)$ represents a time-varying control matrix, which comes from solving the related continuous-time Riccati Equation. The results of the zero controller, LQR controller, and other controllers trained by different methods (as presented in Table~\ref{tab:app_zero_lqr}) are collectively analyzed for a comprehensive comparison. The performance of both zero control and LQR is inferior, indicating the increased difficulty of the respective problems. Notably, the satellite control problem emerges as the least challenging among the scenarios explored in this study, while the quadrotor problem, especially with a longer time horizon, presents a greater level of difficulty.}

% \begin{table}[!ht]
% \centering
% \caption{Statistics of cost ratio of zero control and LQR}
% \label{tab:deter_sat}
% \begin{tabular}{ccccccc}
% \hline
% \textbf{System}              & \textbf{Mean} & \textbf{Std} & \textbf{Max} & \textbf{Min} & \textbf{Median} \\ \hline
% \textbf{Satellite, T=20, $u=0$}          &  118.48             &  86.02            &  407.23            &  12.75                             &   90.65             \\
% \textbf{Satellite, T=20, LQR}          &  1.40             &  0.31            &  3.03            &  1.09                             &   1.33             \\
% \textbf{Quadrotor, T=4, $u=0$}         & 2246.29         & 928.00       & 5179.08      & 1059.37         & 2011.27           \\
% \textbf{Quadrotor, T=4, LQR}         & 93.72         & 35.18       & 201.64      & 48.90         & 84.86           \\
% \textbf{Quadrotor, T=8, $u=0$}         & 229094.04         & 105191.03       & 565065.73      & 94190.64         & 205459.91           \\
% \textbf{Quadrotor, T=8, LQR}         & 974.49         & 424.82       & 2315.50      & 424.25         & 880.33           \\
% \textbf{Quadrotor, T=16, $u=0$}         & 34735478.38         & 17882135.90       & 98257296.49      & 11788478.26         & 31225816.37           \\
% \textbf{Quadrotor, T=16, LQR}         & 9570.70         & 4769.68       & 26384.28      & 3389.41         & 8718.36           \\
% \hline
% \end{tabular}
% \end{table}

\begin{table}[!ht]
\centering
\caption{\revisec{Mean of cost ratio of zero control, finite horizon LQR, direct policy optimization, offline supervised learning and fine-tuned controllers. The last three columns share the same result as Table~\ref{tab:deter_sat} where $\boldsymbol{x}_0 \in {\mathcal{S}}_{\text{sate}}$ and Table~\ref{tab:varying_time} where $\boldsymbol{x}_0 \in \tilde{\mathcal{S}}_{\text{quad}}$.}}
\label{tab:app_zero_lqr}
\begin{tabular}{cccccc}
\hline
\textbf{System}              & \textbf{Zero Control} & \textbf{LQR} & \textbf{DO}  &  \textbf{SL} & \textbf{Fine-tune}\\ \hline
\textbf{Satellite, $T=20$}          &  118.48             &  1.03  &  1.05  & 1.00  & 1.00              \\
\textbf{Quadrotor, $T=4$}         & 2246.29         & 84.53   &  1.05  &  1.01  &  1.00   \\
\textbf{Quadrotor, $T=8$}         & 229094.04         & 1602.78  & 1.63 & 1.15  & 1.03  \\
\textbf{Quadrotor, $T=16$}         & 34735478.38         & 74110.31 & 151.71  & 9.59  &  2.31 \\
\hline
\end{tabular}
\end{table}

\section{Experimental Details} 
All the experiments are run on a single NVIDIA 3090 GPU.

\subsection{Training Details}\label{app:exp_details}

Hyper-parameters are summarized in Table~\ref{tab:hyperparameter}, and we briefly explain some details. 

In our RL experiments, we use a discretization timestep of $\Delta t=0.005$ to build the environment, resulting in 4000 steps per episode. The reward at each timestep $t$ is the running cost integrated over a small time interval, $R_{t} = -\int_t^{t+\Delta t} {L}(\boldsymbol{x}, \boldsymbol{u}) dt$, and the terminal reward additionally subtracts the terminal cost $M(\boldsymbol{x}(T))$. In RL, the goal is to learn the policy $\boldsymbol{u}$ that maximizes the expected discounted return $G_t = \sum_{k=0}^{\infty} \gamma^k R_{t+k+1}$, where the discount factor $\gamma=1.0$ in our finite horizon problem. Note that the expected discounted return in RL to maximize is the total cost in OCP \eqref{eq1-objective} to minimize. 
We use the PPO algorithm implemented in stable-baselines3~\citep{stable-baselines3}.\footnote{https://github.com/DLR-RM/stable-baselines3} This algorithm employs the generalized advantage estimation (GAE, \citealp{gae}) technique to reduce variance in the gradient estimate. We experiment with \texttt{gae\_lambda} from \texttt{\{1.0, 0.99, 0.98, 0.97, 0.96, 0.95\}}, ultimately settling on \texttt{gae\_lambda=0.98} as it yields the best results. To ensure fair comparisons, we set the number of inner loop epochs in PPO to 1, so that the collected data is only used once in RL, just as in direct optimization with dynamics. It is worth noting that the number of trajectories used in direct optimization with dynamics is the product of its batch size and number of iterations (\texttt{1024 $\times$ 2000}), which is the same as the number of episodes used in RL (1x) (\texttt{2048000}), as shown in Table~\ref{tab:hyperparameter}.

We use \texttt{Dopri5} (Dormand--Prince Runge--Kutta method of order (4)5) in the satellite problem and \texttt{RK23} in the quadrotor problem as the ODE solvers. The implementation of direct policy optimization in this paper includes the aforementioned advanced techniques of adaptive solvers and adding checkpoints to achieve desired performance. The adaptive solvers in training are implemented by \cite{ACA_NODE}.\footnote{https://github.com/juntang-zhuang/torch-ACA} 
Since our main results are architecture-agnostic, we use a fully connected neural network as the closed-loop controller, which takes both time and state as the input, \textit{i.e.}, $\boldsymbol{u}^{\mathrm{NN}}(t, \boldsymbol{x})$. 
We use time-marching in satellite and space-marching in quadrotor to successfully generate the BVP solution. We randomly generate a validation dataset of 100 initial states for each problem, and adopt BVP solvers to compute their optimal costs for comparison. As for IVP-enhanced sampling, we follow the same settings as \cite{IVP_Enhanced} does, with temporal grid points on $0<10<14<16$. The number of initial states in uniform sampling and adaptive sampling is the same. %
We remark that a learning rate decay strategy is incorporated to achieve good performance in SL, following common practices in the field. However, in the case of direct policy optimization, we are unable to employ similar learning rate decay due to its inefficient optimization process. Given the limited number of iterations currently used, reducing the learning rate in direct policy optimization would significantly slow down the training process and 
hurt performance. On the other hand, at the fine-tuning stage, we have the flexibility to utilize much smaller learning rates, allowing us to achieve optimal performance.

We present further detailed information on the analysis of optimization landscapes in Figure~\ref{fig:opt_landscape}, following the settings outlined in \cite{bnopt}. Our focus is primarily on local evaluations conducted during the training process. We denote the loss function as $l(\theta)$ and the learning rate (used in the supervised learning or direct policy optimization) as $lr$. During the local evaluation, we perform updates similar to stochastic gradient descent (SGD), $\theta^\prime = \hat{\theta} - {step\_size} \times \nabla_{\theta} l(\hat{\theta})$, where $step\_size = step\_ratio \times lr$ and $\hat{\theta}$ denotes the current parameter in training. Specifically, the multiple values of $step\_ratio$ selected for supervised learning (SL) are within [1/100, 100], while for direct optimization, the values of $step\_ratio$ are within [1/100, 1/10]. Notably, in direct optimization \reviseb{which is randomly initialized}, we are limited to choosing $step\_ratio$ values much smaller than 1.0. This restriction arises due to the differing update rules between Adam (the training optimizer) and SGD (the update rule in local evaluation). The need for this careful adjustment is driven by the high sensitivity of direct policy optimization to the choice of optimizer and learning rate, particularly in the initial stages when gradients are significantly larger. \reviseb{The $step\_ratio$ is not strictly limited to such small values when the weights are pre-trained and have a much smoother landscape. We remark that we use the same $step\_ratio$ in direct optimization for both randomly initialized and with pre-trained weights.} By selecting appropriate $step\_ratio$ values within this constrained range, we ensure a valid evaluation process.

\begin{table}[!htbp]
\small
\caption{Experimental Settings}
\label{tab:hyperparameter}
\begin{tabular}{ll}
\hline
\textbf{Hyper-parameter}                   & \textbf{Setting}                                            \\ \hline
Network architecture                       & Three hidden layers with a hidden size of 64 \\
Activation function                       & \texttt{Tanh} \\
ODE solver                               &  \texttt{Dopri5} with \texttt{atol=1e-5} and \texttt{rtol=1e-5} (Satellite)                        \\
                                           & \texttt{RK23} with \texttt{atol=1e-5} and \texttt{rtol=1e-5} (Quadrotor)                             \\
BVP technique               & Time-marching~\citep{Adaptive_BVP_HJB} (Satellite)                                             \\
                                           & Space-marching~\citep{ML_Enhanced_Landing} (Quadrotor)                                \\
Trajectories  for evaluation          & 100                                                         \\
Trajectories for training             & 100 (Satellite)                                             \\
                                           & 500 (Quadrotor in $\tilde{\mathcal{S}}_{quad}$)                                \\
                                           & 1000 (Quadrotor in $\mathcal{S}_{quad}$)                               \\
Epochs in SL    & 100 (Satellite)                                             \\
                                           & 1000 (Quadrotor in $\tilde{\mathcal{S}}_{quad}$)                                        \\
                                           & 2000 (Quadrotor in $\mathcal{S}_{quad}$)                                        \\
Iterations in DO & 2000 (Satellite)                                            \\
                                           & 3000 (Quadrotor)                                            \\
Iterations in   Fine-tuning         & 100 (Quadrotor in $\tilde{\mathcal{S}}_{quad}$)                                \\
                                           & 100 / 1000 (Quadrotor in $\mathcal{S}_{quad}$)                         \\
Episodes in RL              &  2048000 / 10240000                                                 \\
Optimizer                                  & Adam~\citep{kingma2015adam}                                                        \\
Learning rate in SL    & 0.01 (Satellite)                                \\
                                           & 0.001 (Quadrotor in $\tilde{\mathcal{S}}_{quad}$)                                \\
                                           & 0.01 and decays every 500 epochs by 0.5 (uniform in $\mathcal{S}_{quad}$)        \\
                                           & 0.005 and decays every 500 epochs by 0.5 (adaptive in $\mathcal{S}_{quad}$)       \\
Learning rate in DO     & 0.01                                                        \\
Learning rate in Fine-tuning               & 0.0001                                                      \\
Learning rate in RL               & 0.0005                                                      \\
Batch size in SL        & 1024 (Satellite)                                                        \\
        & 4096 (Quadrotor)                                                       \\
Batch size in DO        & 1024  (Satellite)                                                       \\
        & 2048 (Quadrotor)                                                        \\
Batch size in Fine-tuning                  & 2048                                                        \\
Batch size in RL                 & 2048                                                        \\ \hline
\end{tabular}
\end{table}

\subsection{Training Time} \label{app:training_time}

The column \texttt{Data} denotes the time for generating the whole dataset by BVP solvers according to the related marching methods. We apply parallel techniques to fasten the total procedure of data generation in quadrotor with 24 processors. The addition time in \texttt{Quad\_Large\_Ada} is related to the network training in adaptive sampling. The other columns represent the training time of supervised learning, the fine-tuning time based on the pre-trained network, and the training time by direct optimization, respectively. The rows \texttt{Quad\_Small\_T4}, \texttt{Quad\_Small\_T8} and \texttt{Quad\_Small\_T16} denote the experiments initiated in the small domain $\tilde{\mathcal{S}}_{\text{quad}}$ with varying time spans of $T=4, 8, 16$. The last two rows are experiments in $\mathcal{S}_{\text{quad}}$ with $T=16$, where \texttt{Quad\_Large\_Uni} denotes that the dataset for supervised learning is generated based on uniform initial states, whilst \texttt{Quad\_Large\_Ada} denotes that the dataset is generated by adaptive sampling.

Table~\ref{tab:training_time} demonstrates that supervised learning costs less time than direct optimization and highlights the priority of the proposed Pre-train and Fine-tune strategy in terms of time consumption. Direct optimization is harder to train and costs more time, especially when the problem is more complex and the time horizon is longer. There is no need for fine-tuning in satellite, since the pre-trained network is already very close to the optimal control. The fine-tuning time reported is only for 100 iterations, which are all within 10 minutes. As shown in Table~\ref{tab:varying_time} and Table~\ref{tab:exp_10_16}, with only 100 iterations of fine-tuning, the performance of the controller can be improved significantly. We also record the fine-tuning time of 1000 iterations, which is about 92 minutes for both the uniform dataset and the adaptive dataset. All experiments related to fine-tuning are fine-tuned for 100 iterations, except the ones mentioned in Table~\ref{tab:exp_10_16}. 

We utilize 256 parallel environments to collect data for reinforcement learning. The training process takes 85 minutes in RL (1x) and 434 minutes in RL (5x). The training times are comparable between RL (1x) and BP with fully-known dynamics when using the same number of trajectories, and are approximately five times longer for RL (5x) at a larger sample scale.

\begin{table}[!htbp]
\caption{Summary on time (min) of data generation and network training}
\label{tab:training_time}
\centering
\begin{tabular}{lllll}
\hline
\textbf{Problem}   & \textbf{Data} & \textbf{Supervised} & \textbf{Fine-tuning} & \textbf{Direct} \\ \hline
\textbf{Satellite}          & 0.8                                         & 0.7                                                 &    ---                                                  & 87                                             \\
\textbf{Quad\_Small\_T4}  &   2.4                                       &  13.5                                                   &   4.4                                               &   106                                        \\
\textbf{Quad\_Small\_T8}  &    2.2                                    &   13.5                                                  &    6.0                                             &     157                                      \\
\textbf{Quad\_Small\_T16}  &   2.2                                       &   13.5                                                    &      7.4                                            &                                           387   \\
\textbf{Quad\_Large\_Uni} &  12.8                                     & 55.5                                                    & 9.1                                                 & 426                                             \\
\textbf{Quad\_Large\_Ada} &   12.8 + 162                                 &       51.5                                              & 9.4                                                 & 426                                             \\ \hline
\end{tabular}
\end{table}

\section{Full Dynamics of Satellite and Quadrotor}\label{app:dyn}

\subsection{Satellite}  \label{app:dyn_sat}
We introduce the full dynamics of the satellite problem~\citep{Adaptive_BVP_HJB, Satellite_15, Satellite_17} that are considered in Section~\ref{sec: exp_sat}.

The state is $\boldsymbol{x}=\left(\boldsymbol{v}^{\mathrm{T}}, \boldsymbol{\omega}^{\mathrm{T}} \right)^{\mathrm{T}}$, where $\boldsymbol{v} = (\phi, \theta, \psi)^{\mathrm{T}}$ denotes the attitude of the satellite in Euler angles and $\boldsymbol{\omega}=\left(\omega_1, \omega_2, \omega_3\right)^{\mathrm{T}}$ represents the angular velocity in the body frame. $\phi$ (roll), $\theta$ (pitch), and $\psi$ (yaw) are the rotation angles around the body frame.

The closed-loop control is $\boldsymbol{u}(t, \boldsymbol{x}):\left[0, T\right] \times \mathbb{R}^6 \rightarrow \mathbb{R}^m$. The output dimension $m$ represents the number of momentum wheels, which is set to 3 as the fully actuated case. Then we define the constant matrix $\boldsymbol{B} \in \mathbb{R}^{3 \times m}$, $\boldsymbol{J} \in \mathbb{R}^{3 \times 3}$ to denote the combination of the inertia matrices of the momentum wheels, and $\boldsymbol{h} \in \mathbb{R}^3$ to denote the rigid body without wheels as

\begin{equation*}
\boldsymbol{B}=\left[\begin{array}{ccc}
1 & 1 / 20 & 1 / 10 \\
1 / 15 & 1 & 1 / 10 \\
1 / 10 & 1 / 15 & 1
\end{array}\right], \quad \boldsymbol{J}=\left[\begin{array}{lll}
2 & 0 & 0 \\
0 & 3 & 0 \\
0 & 0 & 4
\end{array}\right], \quad \boldsymbol{h}=\left[\begin{array}{l}
1 \\
1 \\
1
\end{array}\right].
\end{equation*}

Finally, the full dynamics are
\begin{equation*}
\left[\begin{array}{c}
\dot{\boldsymbol{v}} \\
\boldsymbol{J} \dot{\boldsymbol{\omega}}
\end{array}\right]=\left[\begin{array}{c}
\boldsymbol{E}(\boldsymbol{v}) \boldsymbol{\omega} \\
\boldsymbol{S}(\boldsymbol{\omega}) \boldsymbol{R}(\boldsymbol{v}) \boldsymbol{h}+\boldsymbol{B u}
\end{array}\right],
\end{equation*}

where the skew-symmetric matrix and rotation matrices $\boldsymbol{S}(\boldsymbol{\omega}), \boldsymbol{E}(\boldsymbol{v}),  \boldsymbol{R}(\boldsymbol{v}) : \mathbb{R}^3 \rightarrow \mathbb{R}^{3 \times 3}$ are defined as 
\begin{equation*}
\boldsymbol{S}(\boldsymbol{\omega})=\left[\begin{array}{ccc}
0 & \omega_3 & -\omega_2 \\
-\omega_3 & 0 & \omega_1 \\
\omega_2 & -\omega_1 & 0
\end{array}\right], \quad \boldsymbol{E}(\boldsymbol{v})=\left[\begin{array}{ccc}
1 & \sin \phi \tan \theta & \cos \phi \tan \theta \\
0 & \cos \phi & -\sin \phi \\
0 & \sin \phi / \cos \theta & \cos \phi / \cos \theta
\end{array}\right],
\end{equation*}

\begin{equation*}
\boldsymbol{R}(\boldsymbol{v})=\left[\begin{array}{ccc}
\cos \theta \cos \psi & \cos \theta \sin \psi & -\sin \theta \\
\sin \phi \sin \theta \cos \psi-\cos \phi \sin \psi & \sin \phi \sin \theta \sin \psi+\cos \phi \cos \psi & \cos \theta \sin \phi \\
\cos \phi \sin \theta \cos \psi+\sin \phi \sin \psi & \cos \phi \sin \theta \sin \psi-\sin \phi \cos \psi & \cos \theta \cos \phi
\end{array}\right].
\end{equation*}

The running cost and the final cost are defined as  
\begin{equation*}
L(\boldsymbol{x}, \boldsymbol{u})= W_1 \|\boldsymbol{v}\|^2+ W_2 \|\boldsymbol{\omega}\|^2+ W_3 \|\boldsymbol{u}\|^2,
\end{equation*}
\begin{equation*}
M(\boldsymbol{x}(T)) = {W_4}\left\|\boldsymbol{v}(T)\right\|^2+ {W_5}\left\|\boldsymbol{\omega}(T)\right\|^2,
\end{equation*}
where the weights are set as $W_1=\frac{1}{2}, W_2=5, W_3=\frac{1}{4}, W_4=\frac{1}{2}, W_5=\frac{1}{2}$.

\subsection{Quadrotor} \label{app:dyn_qua}
We describe the full dynamics of the optimal landing problem of a quadrotor~\citep{ML_Enhanced_Landing, IVP_Enhanced, Quadrotor_04, Quadrotor_06, Quadrotor_12} which are considered in Section~\ref{sec:exp_quadrotor}.

The state is $\boldsymbol{x} = \left(\boldsymbol{p}^{\mathrm{T}}, \boldsymbol{v}_b^{\mathrm{T}}, \boldsymbol{\eta}^{\mathrm{T}}, \boldsymbol{w}_b^{\mathrm{T}}\right)^{\mathrm{T}}\in \mathbb{R}^{12}$ where $\boldsymbol{p}=(x, y, z)^{\mathrm{T}} \in \mathbb{R}^3$ denotes the position in Earth-fixed coordinates, $\boldsymbol{v}_b=(v_x, v_y, v_z)^{\mathrm{T}} \in \mathbb{R}^3$ denotes the velocity with respect to the body frame, $\boldsymbol{\eta}=(\phi, \theta, \psi)^{\mathrm{T}} \in \mathbb{R}^3$ denotes the attitude in Earth-fixed coordinates, and $\boldsymbol{w}_b \in \mathbb{R}^3$ denotes the angular velocity in the body frame.

The closed-loop control is $\boldsymbol{u}(t, \boldsymbol{x}):\left[0, T\right] \times \mathbb{R}^{12} \rightarrow \mathbb{R}^4$. In practice, the individual rotor thrusts $\boldsymbol{F} = \left(F_1, F_2, F_3, F_4\right)^{\mathrm{T}}$ is applied to steer the quadrotor, where $\boldsymbol{u}=E \boldsymbol{F}$. We use $l$ as the distance from the rotor to the quadrotor's center of gravity, and $c$ as a constant that relates the rotor angular momentum to the rotor thrust (normal force), then the matrix $E$ is defined as

\begin{equation*}
E=\left[\begin{array}{cccc}
1 & 1 & 1 & 1 \\
0 & l & 0 & -l \\
-l & 0 & l & 0 \\
c & -c & c & -c
\end{array}\right].
\end{equation*}
With $\boldsymbol{F}^*=E^{-1} \boldsymbol{u}^*$, the optimal $\boldsymbol{F}^*$ can be directly computed, once we obtain the optimal control $\boldsymbol{u}^*$. 

Using the same parameters as~\cite{Quadrotor_06}, the mass and the inertia matrix are defined as $m=2 k g$ and $\boldsymbol{J}=\operatorname{diag}\left(J_x, J_y, J_z\right)$, where $J_x=J_y=\frac{1}{2} J_z=1.2416 \mathrm{~kg} \cdot \mathrm{m}^2$. Another constant vector is $\boldsymbol{g}=(0, 0, g)^{\mathrm{T}}$, where $g=9.81 \mathrm{~m} / \mathrm{s}^2$ denotes the acceleration of gravity on Earth. The rotation matrix $\boldsymbol{R}(\boldsymbol{\eta}) \in S O(3)$ denotes the transformation from the Earth-fixed coordinates to the body-fixed coordinates, and the attitude kinematic matrix $\boldsymbol{K}(\boldsymbol{\eta})$ relates the time derivate of the attitude with the associated angular rate as

\begin{equation*}
\boldsymbol{R}(\boldsymbol{\eta})=\left[\begin{array}{ccc}
\cos \theta \cos \psi & \cos \theta \sin \psi & -\sin \theta \\
\sin \theta \cos \psi \sin \phi-\sin \psi \cos \phi & \sin \theta \sin \psi \sin \phi+\cos \psi \cos \phi & \cos \theta \sin \phi \\
\sin \theta \cos \psi \cos \phi+\sin \psi \sin \phi & \sin \theta \sin \psi \cos \phi-\cos \psi \sin \phi & \cos \theta \cos \phi
\end{array}\right],
\end{equation*}

\begin{equation*}
\boldsymbol{K}(\boldsymbol{\eta})=\left[\begin{array}{ccc}
1 & \sin \phi \tan \theta & \cos \phi \tan \theta \\
0 & \cos \phi & -\sin \phi \\
0 & \sin \phi \sec \theta & \cos \phi \sec \theta
\end{array}\right].
\end{equation*}

Finally, the dynamics are

\begin{equation*}
\left\{\begin{array}{l}
\dot{\boldsymbol{p}}=\boldsymbol{R}^\mathrm{T}(\boldsymbol{\eta}) \boldsymbol{v}_b \\
\dot{\boldsymbol{v}}_b=-\boldsymbol{w}_b \times \boldsymbol{v}_b-\boldsymbol{R}(\boldsymbol{\eta}) \boldsymbol{g}+\frac{1}{m} A \boldsymbol{u} \\
\dot{\boldsymbol{\eta}}=\boldsymbol{K}(\boldsymbol{\eta}) \boldsymbol{w}_b \\
\dot{\boldsymbol{w}}_b=-\boldsymbol{J}^{-1} \boldsymbol{w}_b \times \boldsymbol{J} \boldsymbol{w}_b+\boldsymbol{J}^{-1} B \boldsymbol{u}
\end{array}\right. ,
\end{equation*}
where the constant matrices $A$ and $B$ are defined as
\begin{equation*}
A=\left[\begin{array}{llll}
0 & 0 & 0 & 0 \\
0 & 0 & 0 & 0 \\
1 & 0 & 0 & 0
\end{array}\right], \quad B=\left[\begin{array}{llll}
0 & 1 & 0 & 0 \\
0 & 0 & 1 & 0 \\
0 & 0 & 0 & 1
\end{array}\right].
\end{equation*}

The running cost and the final cost are defined as  
\begin{equation*}
L(\boldsymbol{x}, \boldsymbol{u})=\left(\boldsymbol{u}-\boldsymbol{u}_d\right)^\mathrm{T} Q_u\left(\boldsymbol{u}-\boldsymbol{u}_d\right),
\end{equation*}
\begin{equation*}
M(\boldsymbol{x})=\boldsymbol{p}^\mathrm{T} Q_{p f} \boldsymbol{p}+\boldsymbol{v}^\mathrm{T} Q_{v f} \boldsymbol{v}+\boldsymbol{\eta}^\mathrm{T} Q_{\eta f} \boldsymbol{\eta}+\boldsymbol{w}^\mathrm{T} Q_{w f} \boldsymbol{w}=\boldsymbol{x}^\mathrm{T} Q_f \boldsymbol{x},
\end{equation*}
where $Q_u=\operatorname{diag}(1,1,1,1)$, $Q_{p f}=5 I_3$, $Q_{v f}=10 I_3$, $Q_{\eta f}=25 I_3$ and $Q_{w f}=50 I_3$, where $I_3$ denotes the identity matrix in three-dimensional space.

\bibliographystyle{elsarticle-num-names} 
\bibliography{ref}

\begin{thebibliography}{47}
\expandafter\ifx\csname natexlab\endcsname\relax\def\natexlab#1{#1}\fi
\providecommand{\url}[1]{\texttt{#1}}
\providecommand{\href}[2]{#2}
\providecommand{\path}[1]{#1}
\providecommand{\DOIprefix}{doi:}
\providecommand{\ArXivprefix}{arXiv:}
\providecommand{\URLprefix}{URL: }
\providecommand{\Pubmedprefix}{pmid:}
\providecommand{\doi}[1]{\href{http://dx.doi.org/#1}{\path{#1}}}
\providecommand{\Pubmed}[1]{\href{pmid:#1}{\path{#1}}}
\providecommand{\bibinfo}[2]{#2}
\ifx\xfnm\relax \def\xfnm[#1]{\unskip,\space#1}\fi
%Type = Book
\bibitem[{Franklin et~al.(2020)Franklin, Powell, and
  Emami-Naeini}]{franklinFeedbackControlDynamic2020}
\bibinfo{author}{G.~F. Franklin}, \bibinfo{author}{J.~D. Powell},
  \bibinfo{author}{A.~Emami-Naeini}, \bibinfo{title}{Feedback Control of
  Dynamic Systems}, \bibinfo{edition}{eighth} ed., \bibinfo{publisher}{{Pearson
  Education Limited}}, \bibinfo{year}{2020}.
%Type = Article
\bibitem[{Han and E(2016)}]{hanDeepLearningApproximation2016}
\bibinfo{author}{J.~Han}, \bibinfo{author}{W.~E},
\newblock \bibinfo{title}{Deep learning approximation for stochastic control
  problems},
\newblock \bibinfo{journal}{arXiv preprint arXiv:1611.07422}
  (\bibinfo{year}{2016}).
%Type = Article
\bibitem[{Nakamura-Zimmerer et~al.(2021)Nakamura-Zimmerer, Gong, and
  Kang}]{Adaptive_BVP_HJB}
\bibinfo{author}{T.~Nakamura-Zimmerer}, \bibinfo{author}{Q.~Gong},
  \bibinfo{author}{W.~Kang},
\newblock \bibinfo{title}{Adaptive deep learning for high-dimensional
  {Hamilton--Jacobi--Bellman} equations},
\newblock \bibinfo{journal}{SIAM Journal on Scientific Computing}
  \bibinfo{volume}{43} (\bibinfo{year}{2021}) \bibinfo{pages}{A1221--A1247}.
%Type = Article
\bibitem[{Böttcher et~al.(2022)Böttcher, Antulov-Fantulin, and
  Asikis}]{AIPontryagin}
\bibinfo{author}{L.~Böttcher}, \bibinfo{author}{N.~Antulov-Fantulin},
  \bibinfo{author}{T.~Asikis},
\newblock \bibinfo{title}{{{AI Pontryagin}} or how artificial neural networks
  learn to control dynamical systems},
\newblock \bibinfo{journal}{Nature Communications} \bibinfo{volume}{13}
  (\bibinfo{year}{2022}) \bibinfo{pages}{333}.
%Type = Article
\bibitem[{E et~al.(2022)E, Han, and Long}]{Empowering}
\bibinfo{author}{W.~E}, \bibinfo{author}{J.~Han}, \bibinfo{author}{J.~Long},
\newblock \bibinfo{title}{Empowering optimal control with machine learning: A
  perspective from model predictive control},
\newblock \bibinfo{journal}{arXiv preprint arXiv:2205.07990}
  (\bibinfo{year}{2022}).
%Type = Inproceedings
\bibitem[{Ainsworth et~al.(2021)Ainsworth, Lowrey, Thickstun, Harchaoui, and
  Srinivasa}]{CTPG}
\bibinfo{author}{S.~Ainsworth}, \bibinfo{author}{K.~Lowrey},
  \bibinfo{author}{J.~Thickstun}, \bibinfo{author}{Z.~Harchaoui},
  \bibinfo{author}{S.~Srinivasa},
\newblock \bibinfo{title}{Faster policy learning with continuous-time
  gradients},
\newblock in: \bibinfo{booktitle}{Learning for Dynamics and Control},
  \bibinfo{organization}{PMLR}, \bibinfo{year}{2021}, pp.
  \bibinfo{pages}{1054--1067}.
%Type = Book
\bibitem[{Liberzon(2011)}]{Openloop}
\bibinfo{author}{D.~Liberzon}, \bibinfo{title}{Calculus of variations and
  optimal control theory: a concise introduction},
  \bibinfo{publisher}{Princeton University Press}, \bibinfo{year}{2011}.
%Type = Book
\bibitem[{Bellman(1957)}]{bellmanDynamicProgramming1957}
\bibinfo{author}{R.~Bellman}, \bibinfo{title}{Dynamic Programming}, Rand
  {{Coperation}} Research Study, \bibinfo{publisher}{{Princeton University
  Press}}, \bibinfo{year}{1957}.
%Type = Article
\bibitem[{Azmi et~al.(2021)Azmi, Kalise, and Kunisch}]{azmi2021optimal}
\bibinfo{author}{B.~Azmi}, \bibinfo{author}{D.~Kalise},
  \bibinfo{author}{K.~Kunisch},
\newblock \bibinfo{title}{Optimal feedback law recovery by gradient-augmented
  sparse polynomial regression},
\newblock \bibinfo{journal}{The Journal of Machine Learning Research}
  \bibinfo{volume}{22} (\bibinfo{year}{2021}) \bibinfo{pages}{2205--2236}.
%Type = Article
\bibitem[{Kunisch et~al.(2023)Kunisch, V{\'a}squez-Varas, and
  Walter}]{kunisch2023learning}
\bibinfo{author}{K.~Kunisch}, \bibinfo{author}{D.~V{\'a}squez-Varas},
  \bibinfo{author}{D.~Walter},
\newblock \bibinfo{title}{Learning optimal feedback operators and their sparse
  polynomial approximations},
\newblock \bibinfo{journal}{Journal of Machine Learning Research}
  \bibinfo{volume}{24} (\bibinfo{year}{2023}) \bibinfo{pages}{1--38}.
%Type = Article
\bibitem[{Kang and Wilcox(2017)}]{Satellite_17}
\bibinfo{author}{W.~Kang}, \bibinfo{author}{L.~C. Wilcox},
\newblock \bibinfo{title}{Mitigating the curse of dimensionality: sparse grid
  characteristics method for optimal feedback control and {HJB} equations},
\newblock \bibinfo{journal}{Computational Optimization and Applications}
  \bibinfo{volume}{68} (\bibinfo{year}{2017}) \bibinfo{pages}{289--315}.
%Type = Article
\bibitem[{Weston et~al.(2002)Weston, Chapelle, Vapnik, Elisseeff, and
  Sch{\"o}lkopf}]{weston2002kernel}
\bibinfo{author}{J.~Weston}, \bibinfo{author}{O.~Chapelle},
  \bibinfo{author}{V.~Vapnik}, \bibinfo{author}{A.~Elisseeff},
  \bibinfo{author}{B.~Sch{\"o}lkopf},
\newblock \bibinfo{title}{Kernel dependency estimation},
\newblock \bibinfo{journal}{Advances in neural information processing systems}
  \bibinfo{volume}{15} (\bibinfo{year}{2002}).
%Type = Article
\bibitem[{Meng et~al.(2022)Meng, Zhang, Darbon, and
  Karniadakis}]{meng2022sympocnet}
\bibinfo{author}{T.~Meng}, \bibinfo{author}{Z.~Zhang},
  \bibinfo{author}{J.~Darbon}, \bibinfo{author}{G.~Karniadakis},
\newblock \bibinfo{title}{Sympocnet: Solving optimal control problems with
  applications to high-dimensional multiagent path planning problems},
\newblock \bibinfo{journal}{SIAM Journal on Scientific Computing}
  \bibinfo{volume}{44} (\bibinfo{year}{2022}) \bibinfo{pages}{B1341--B1368}.
%Type = Article
\bibitem[{Onken et~al.(2022)Onken, Nurbekyan, Li, Fung, Osher, and
  Ruthotto}]{onken2022neural}
\bibinfo{author}{D.~Onken}, \bibinfo{author}{L.~Nurbekyan},
  \bibinfo{author}{X.~Li}, \bibinfo{author}{S.~W. Fung},
  \bibinfo{author}{S.~Osher}, \bibinfo{author}{L.~Ruthotto},
\newblock \bibinfo{title}{A neural network approach for high-dimensional
  optimal control applied to multiagent path finding},
\newblock \bibinfo{journal}{IEEE Transactions on Control Systems Technology}
  \bibinfo{volume}{31} (\bibinfo{year}{2022}) \bibinfo{pages}{235--251}.
%Type = Article
\bibitem[{Kierzenka and Shampine(2001)}]{BVP}
\bibinfo{author}{J.~Kierzenka}, \bibinfo{author}{L.~F. Shampine},
\newblock \bibinfo{title}{A bvp solver based on residual control and the maltab
  pse},
\newblock \bibinfo{journal}{ACM Transactions on Mathematical Software (TOMS)}
  \bibinfo{volume}{27} (\bibinfo{year}{2001}) \bibinfo{pages}{299--316}.
%Type = Book
\bibitem[{Jacobson and Mayne(1970)}]{DDP}
\bibinfo{author}{D.~H. Jacobson}, \bibinfo{author}{D.~Q. Mayne},
  \bibinfo{title}{Differential dynamic programming},
  \bibinfo{publisher}{Elsevier Publishing Company}, \bibinfo{year}{1970}.
%Type = Article
\bibitem[{Osa et~al.(2018)Osa, Pajarinen, Neumann, Bagnell, Abbeel, Peters
  et~al.}]{osa2018algorithmic}
\bibinfo{author}{T.~Osa}, \bibinfo{author}{J.~Pajarinen},
  \bibinfo{author}{G.~Neumann}, \bibinfo{author}{J.~A. Bagnell},
  \bibinfo{author}{P.~Abbeel}, \bibinfo{author}{J.~Peters}, et~al.,
\newblock \bibinfo{title}{An algorithmic perspective on imitation learning},
\newblock \bibinfo{journal}{Foundations and Trends{\textregistered} in
  Robotics} \bibinfo{volume}{7} (\bibinfo{year}{2018}) \bibinfo{pages}{1--179}.
%Type = Inproceedings
\bibitem[{Bain and Sammut(1995)}]{bain1995framework}
\bibinfo{author}{M.~Bain}, \bibinfo{author}{C.~Sammut},
\newblock \bibinfo{title}{A framework for behavioural cloning.},
\newblock in: \bibinfo{booktitle}{Machine Intelligence 15},
  \bibinfo{year}{1995}, pp. \bibinfo{pages}{103--129}.
%Type = Article
\bibitem[{Bock and Plitt(1984)}]{Bock_Direct}
\bibinfo{author}{H.~Bock}, \bibinfo{author}{K.~Plitt},
\newblock \bibinfo{title}{A multiple shooting algorithm for direct solution of
  optimal control problems},
\newblock \bibinfo{journal}{IFAC Proceedings Volumes} \bibinfo{volume}{17}
  (\bibinfo{year}{1984}) \bibinfo{pages}{1603--1608}.
%Type = Article
\bibitem[{Betts(1998)}]{Betts_Direct}
\bibinfo{author}{J.~T. Betts},
\newblock \bibinfo{title}{Survey of numerical methods for trajectory
  optimization},
\newblock \bibinfo{journal}{Journal of Guidance, Control, and Dynamics}
  \bibinfo{volume}{21} (\bibinfo{year}{1998}) \bibinfo{pages}{193--207}.
%Type = Article
\bibitem[{Ross and Fahroo(2002)}]{Ross_Direct}
\bibinfo{author}{I.~M. Ross}, \bibinfo{author}{F.~Fahroo},
\newblock \bibinfo{title}{A direct method for solving nonsmooth optimal control
  problems},
\newblock \bibinfo{journal}{IFAC Proceedings Volumes} \bibinfo{volume}{35}
  (\bibinfo{year}{2002}) \bibinfo{pages}{479--484}.
%Type = Incollection
\bibitem[{Diehl et~al.(2006)Diehl, Bock, Diedam, and Wieber}]{Diehl_Direct}
\bibinfo{author}{M.~Diehl}, \bibinfo{author}{H.~G. Bock},
  \bibinfo{author}{H.~Diedam}, \bibinfo{author}{P.-B. Wieber},
\newblock \bibinfo{title}{Fast direct multiple shooting algorithms for optimal
  robot control},
\newblock in: \bibinfo{booktitle}{Fast motions in biomechanics and robotics},
  \bibinfo{publisher}{Springer}, \bibinfo{year}{2006}, pp.
  \bibinfo{pages}{65--93}.
%Type = Article
\bibitem[{Kunisch and Walter(2021)}]{kunisch2021semiglobal}
\bibinfo{author}{K.~Kunisch}, \bibinfo{author}{D.~Walter},
\newblock \bibinfo{title}{Semiglobal optimal feedback stabilization of
  autonomous systems via deep neural network approximation},
\newblock \bibinfo{journal}{ESAIM: Control, Optimisation and Calculus of
  Variations} \bibinfo{volume}{27} (\bibinfo{year}{2021}) \bibinfo{pages}{16}.
%Type = Article
\bibitem[{Rumelhart et~al.(1986)Rumelhart, Hinton, and Williams}]{bp}
\bibinfo{author}{D.~E. Rumelhart}, \bibinfo{author}{G.~E. Hinton},
  \bibinfo{author}{R.~J. Williams},
\newblock \bibinfo{title}{Learning representations by back-propagating errors},
\newblock \bibinfo{journal}{Nature} \bibinfo{volume}{323}
  (\bibinfo{year}{1986}) \bibinfo{pages}{533--536}.
%Type = Book
\bibitem[{Pontryagin et~al.(1962)Pontryagin, Mishchenko, and
  RV}]{pontryaginMathematicalTheoryOptimal1962}
\bibinfo{author}{L.~S. Pontryagin}, \bibinfo{author}{E.~Mishchenko},
  \bibinfo{author}{G.~RV}, \bibinfo{title}{The mathematical theory of optimal
  processes}, \bibinfo{year}{1962}.
%Type = Article
\bibitem[{Schulman et~al.(2017)Schulman, Wolski, Dhariwal, Radford, and
  Klimov}]{ppo}
\bibinfo{author}{J.~Schulman}, \bibinfo{author}{F.~Wolski},
  \bibinfo{author}{P.~Dhariwal}, \bibinfo{author}{A.~Radford},
  \bibinfo{author}{O.~Klimov},
\newblock \bibinfo{title}{Proximal policy optimization algorithms},
\newblock \bibinfo{journal}{arXiv preprint arXiv:1707.06347}
  (\bibinfo{year}{2017}).
%Type = Inproceedings
\bibitem[{Zang et~al.(2022)Zang, Long, Zhang, Hu, E, and
  Han}]{ML_Enhanced_Landing}
\bibinfo{author}{Y.~Zang}, \bibinfo{author}{J.~Long},
  \bibinfo{author}{X.~Zhang}, \bibinfo{author}{W.~Hu}, \bibinfo{author}{W.~E},
  \bibinfo{author}{J.~Han},
\newblock \bibinfo{title}{A machine learning enhanced algorithm for the optimal
  landing problem},
\newblock in: \bibinfo{booktitle}{3rd Annual Conference on Mathematical and
  Scientific Machine Learning}, \bibinfo{organization}{PMLR},
  \bibinfo{year}{2022}, pp. \bibinfo{pages}{1--20}.
%Type = Article
\bibitem[{Long and Han(2022)}]{long2022perturbational}
\bibinfo{author}{J.~Long}, \bibinfo{author}{J.~Han},
\newblock \bibinfo{title}{Perturbational complexity by distribution mismatch: A
  systematic analysis of reinforcement learning in reproducing kernel hilbert
  space},
\newblock \bibinfo{journal}{Journal of Machine Learning vol}
  \bibinfo{volume}{1} (\bibinfo{year}{2022}) \bibinfo{pages}{1--34}.
%Type = Article
\bibitem[{Zhang et~al.(2022)Zhang, Long, Hu, E, and Han}]{IVP_Enhanced}
\bibinfo{author}{X.~Zhang}, \bibinfo{author}{J.~Long}, \bibinfo{author}{W.~Hu},
  \bibinfo{author}{W.~E}, \bibinfo{author}{J.~Han},
\newblock \bibinfo{title}{Initial value problem enhanced sampling for
  closed-loop optimal control design with deep neural networks},
\newblock \bibinfo{journal}{arXiv preprint arXiv:2209.04078}
  (\bibinfo{year}{2022}).
%Type = Article
\bibitem[{Nair et~al.(2020)Nair, Gupta, Dalal, and Levine}]{nair2020offline}
\bibinfo{author}{A.~Nair}, \bibinfo{author}{A.~Gupta},
  \bibinfo{author}{M.~Dalal}, \bibinfo{author}{S.~Levine},
\newblock \bibinfo{title}{Awac: Accelerating online reinforcement learning with
  offline datasets},
\newblock \bibinfo{journal}{arXiv preprint arXiv:2006.09359}
  (\bibinfo{year}{2020}).
%Type = Inproceedings
\bibitem[{Lee et~al.(2022)Lee, Seo, Lee, Abbeel, and Shin}]{lee2022offline}
\bibinfo{author}{S.~Lee}, \bibinfo{author}{Y.~Seo}, \bibinfo{author}{K.~Lee},
  \bibinfo{author}{P.~Abbeel}, \bibinfo{author}{J.~Shin},
\newblock \bibinfo{title}{Offline-to-online reinforcement learning via balanced
  replay and pessimistic q-ensemble},
\newblock in: \bibinfo{booktitle}{Conference on Robot Learning},
  \bibinfo{organization}{PMLR}, \bibinfo{year}{2022}, pp.
  \bibinfo{pages}{1702--1712}.
%Type = Article
\bibitem[{Levine et~al.(2020)Levine, Kumar, Tucker, and Fu}]{levine2020offline}
\bibinfo{author}{S.~Levine}, \bibinfo{author}{A.~Kumar},
  \bibinfo{author}{G.~Tucker}, \bibinfo{author}{J.~Fu},
\newblock \bibinfo{title}{Offline reinforcement learning: Tutorial, review, and
  perspectives on open problems},
\newblock \bibinfo{journal}{arXiv preprint arXiv:2005.01643}
  (\bibinfo{year}{2020}).
%Type = Inproceedings
\bibitem[{Kang and Wilcox(2015)}]{Satellite_15}
\bibinfo{author}{W.~Kang}, \bibinfo{author}{L.~Wilcox},
\newblock \bibinfo{title}{A causality free computational method for {HJB}
  equations with application to rigid body satellites},
\newblock in: \bibinfo{booktitle}{AIAA Guidance, Navigation, and Control
  Conference}, \bibinfo{year}{2015}, p. \bibinfo{pages}{2009}.
%Type = Inproceedings
\bibitem[{Bouabdallah et~al.(2004)Bouabdallah, Murrieri, and
  Siegwart}]{Quadrotor_04}
\bibinfo{author}{S.~Bouabdallah}, \bibinfo{author}{P.~Murrieri},
  \bibinfo{author}{R.~Siegwart},
\newblock \bibinfo{title}{Design and control of an indoor micro quadrotor},
\newblock in: \bibinfo{booktitle}{IEEE International Conference on Robotics and
  Automation}, volume~\bibinfo{volume}{5}, \bibinfo{organization}{IEEE},
  \bibinfo{year}{2004}, pp. \bibinfo{pages}{4393--4398}.
%Type = Inproceedings
\bibitem[{Madani and Benallegue(2006)}]{Quadrotor_06}
\bibinfo{author}{T.~Madani}, \bibinfo{author}{A.~Benallegue},
\newblock \bibinfo{title}{Control of a quadrotor mini-helicopter via full state
  backstepping technique},
\newblock in: \bibinfo{booktitle}{Proceedings of the 45th IEEE Conference on
  Decision and Control}, \bibinfo{organization}{IEEE}, \bibinfo{year}{2006},
  pp. \bibinfo{pages}{1515--1520}.
%Type = Article
\bibitem[{Mahony et~al.(2012)Mahony, Kumar, and Corke}]{Quadrotor_12}
\bibinfo{author}{R.~Mahony}, \bibinfo{author}{V.~Kumar},
  \bibinfo{author}{P.~Corke},
\newblock \bibinfo{title}{Multirotor aerial vehicles: Modeling, estimation, and
  control of quadrotor},
\newblock \bibinfo{journal}{IEEE Robotics and Automation magazine}
  \bibinfo{volume}{19} (\bibinfo{year}{2012}) \bibinfo{pages}{20--32}.
%Type = Inproceedings
\bibitem[{Santurkar et~al.(2018)Santurkar, Tsipras, Ilyas, and Madry}]{bnopt}
\bibinfo{author}{S.~Santurkar}, \bibinfo{author}{D.~Tsipras},
  \bibinfo{author}{A.~Ilyas}, \bibinfo{author}{A.~Madry},
\newblock \bibinfo{title}{How does batch normalization help optimization?},
\newblock in: \bibinfo{booktitle}{Advances in Neural Information Processing
  Systems}, volume~\bibinfo{volume}{32}, \bibinfo{year}{2018}.
%Type = Inproceedings
\bibitem[{Zhuang et~al.(2020)Zhuang, Dvornek, Li, Tatikonda, Papademetris, and
  Duncan}]{ACA_NODE}
\bibinfo{author}{J.~Zhuang}, \bibinfo{author}{N.~Dvornek},
  \bibinfo{author}{X.~Li}, \bibinfo{author}{S.~Tatikonda},
  \bibinfo{author}{X.~Papademetris}, \bibinfo{author}{J.~Duncan},
\newblock \bibinfo{title}{Adaptive checkpoint adjoint method for gradient
  estimation in neural {ODE}},
\newblock in: \bibinfo{booktitle}{International Conference on Machine
  Learning}, \bibinfo{organization}{PMLR}, \bibinfo{year}{2020}, pp.
  \bibinfo{pages}{11639--11649}.
%Type = Inproceedings
\bibitem[{Chen et~al.(2018)Chen, Rubanova, Bettencourt, and
  Duvenaud}]{NeuralODE}
\bibinfo{author}{R.~T.~Q. Chen}, \bibinfo{author}{Y.~Rubanova},
  \bibinfo{author}{J.~Bettencourt}, \bibinfo{author}{D.~Duvenaud},
\newblock \bibinfo{title}{Neural ordinary differential equations},
\newblock in: \bibinfo{booktitle}{Advances in Neural Information Processing
  Systems}, volume~\bibinfo{volume}{32}, \bibinfo{year}{2018}.
%Type = Article
\bibitem[{Gholami et~al.(2019)Gholami, Keutzer, and Biros}]{ANODE}
\bibinfo{author}{A.~Gholami}, \bibinfo{author}{K.~Keutzer},
  \bibinfo{author}{G.~Biros},
\newblock \bibinfo{title}{Anode: Unconditionally accurate memory-efficient
  gradients for neural odes},
\newblock \bibinfo{journal}{arXiv preprint arXiv:1902.10298}
  (\bibinfo{year}{2019}).
%Type = Article
\bibitem[{Miki et~al.(2022)Miki, Lee, Hwangbo, Wellhausen, Koltun, and
  Hutter}]{ppo_robot}
\bibinfo{author}{T.~Miki}, \bibinfo{author}{J.~Lee},
  \bibinfo{author}{J.~Hwangbo}, \bibinfo{author}{L.~Wellhausen},
  \bibinfo{author}{V.~Koltun}, \bibinfo{author}{M.~Hutter},
\newblock \bibinfo{title}{Learning robust perceptive locomotion for quadrupedal
  robots in the wild},
\newblock \bibinfo{journal}{Science Robotics} \bibinfo{volume}{7}
  (\bibinfo{year}{2022}) \bibinfo{pages}{eabk2822}.
%Type = Inproceedings
\bibitem[{Ouyang et~al.(2022)Ouyang, Wu, Jiang, Almeida, Wainwright, Mishkin,
  Zhang, Agarwal, Slama, Ray, Schulman, Hilton, Kelton, Miller, Simens, Askell,
  Welinder, Christiano, Leike, and Lowe}]{instrcutgpt}
\bibinfo{author}{L.~Ouyang}, \bibinfo{author}{J.~Wu},
  \bibinfo{author}{X.~Jiang}, \bibinfo{author}{D.~Almeida},
  \bibinfo{author}{C.~Wainwright}, \bibinfo{author}{P.~Mishkin},
  \bibinfo{author}{C.~Zhang}, \bibinfo{author}{S.~Agarwal},
  \bibinfo{author}{K.~Slama}, \bibinfo{author}{A.~Ray},
  \bibinfo{author}{J.~Schulman}, \bibinfo{author}{J.~Hilton},
  \bibinfo{author}{F.~Kelton}, \bibinfo{author}{L.~Miller},
  \bibinfo{author}{M.~Simens}, \bibinfo{author}{A.~Askell},
  \bibinfo{author}{P.~Welinder}, \bibinfo{author}{P.~F. Christiano},
  \bibinfo{author}{J.~Leike}, \bibinfo{author}{R.~Lowe},
\newblock \bibinfo{title}{Training language models to follow instructions with
  human feedback},
\newblock in: \bibinfo{editor}{S.~Koyejo}, \bibinfo{editor}{S.~Mohamed},
  \bibinfo{editor}{A.~Agarwal}, \bibinfo{editor}{D.~Belgrave},
  \bibinfo{editor}{K.~Cho}, \bibinfo{editor}{A.~Oh} (Eds.),
  \bibinfo{booktitle}{Advances in Neural Information Processing Systems},
  volume~\bibinfo{volume}{35}, \bibinfo{publisher}{Curran Associates, Inc.},
  \bibinfo{year}{2022}, pp. \bibinfo{pages}{27730--27744}.
%Type = Book
\bibitem[{Mehrmann(1991)}]{mehrmann1991autonomous}
\bibinfo{author}{V.~L. Mehrmann}, \bibinfo{title}{The autonomous linear
  quadratic control problem: theory and numerical solution},
  \bibinfo{publisher}{Springer}, \bibinfo{year}{1991}.
%Type = Misc
\bibitem[{Tedrake and the Drake Development~Team(2019)}]{drake}
\bibinfo{author}{R.~Tedrake}, \bibinfo{author}{the Drake Development~Team},
  \bibinfo{title}{Drake: Model-based design and verification for robotics},
  \bibinfo{year}{2019}. \URLprefix \url{https://drake.mit.edu}.
%Type = Article
\bibitem[{Raffin et~al.(2021)Raffin, Hill, Gleave, Kanervisto, Ernestus, and
  Dormann}]{stable-baselines3}
\bibinfo{author}{A.~Raffin}, \bibinfo{author}{A.~Hill},
  \bibinfo{author}{A.~Gleave}, \bibinfo{author}{A.~Kanervisto},
  \bibinfo{author}{M.~Ernestus}, \bibinfo{author}{N.~Dormann},
\newblock \bibinfo{title}{Stable-baselines3: Reliable reinforcement learning
  implementations},
\newblock \bibinfo{journal}{Journal of Machine Learning Research}
  \bibinfo{volume}{22} (\bibinfo{year}{2021}) \bibinfo{pages}{1--8}.
%Type = Inproceedings
\bibitem[{Schulman et~al.(2016)Schulman, Moritz, Levine, Jordan, and
  Abbeel}]{gae}
\bibinfo{author}{J.~Schulman}, \bibinfo{author}{P.~Moritz},
  \bibinfo{author}{S.~Levine}, \bibinfo{author}{M.~Jordan},
  \bibinfo{author}{P.~Abbeel},
\newblock \bibinfo{title}{High-dimensional continuous control using generalized
  advantage estimation},
\newblock in: \bibinfo{booktitle}{Proceedings of the International Conference
  on Learning Representations}, \bibinfo{year}{2016}.
%Type = Inproceedings
\bibitem[{Kingma and Ba(2015)}]{kingma2015adam}
\bibinfo{author}{D.~P. Kingma}, \bibinfo{author}{J.~Ba},
\newblock \bibinfo{title}{Adam: a method for stochastic optimization},
\newblock in: \bibinfo{booktitle}{Proceedings of the International Conference
  on Learning Representations}, \bibinfo{year}{2015}.

\end{thebibliography}

\end{document}